\documentclass[oneside,english]{amsart}
\usepackage[T1]{fontenc}
\usepackage[latin9]{inputenc}
\usepackage{amsthm,amsmath,amssymb,color,graphicx,babel}

\newtheorem{theorem}{Theorem}
\newtheorem{lemma}[theorem]{Lemma}
\newtheorem{corollary}[theorem]{Corollary}

\newtheorem{proposition}[theorem]{Proposition}
\theoremstyle{definition}
\newtheorem{definition}[theorem]{Definition}
\newtheorem{remark}[theorem]{Remark}

\newcommand{\mI}{{\mathcal I}}
\newcommand{\mX}{{\mathcal{S}}}
\newcommand{\mM}{{\mathcal M}}
\newcommand{\mR}{{\mathcal R}}
\newcommand{\mF}{{\mathcal F}}
\newcommand{\mW}{{\mathcal{W}}}
\newcommand{\ZZ}{{\mathbb{Z}}}
\newcommand{\NN}{{\mathbb{N}}}
\newcommand{\mP}{{\mathcal{P}}}
\newcommand{\trho}{{\widetilde{\rho}_0}}

\newcommand{\Cov}{\hbox{\rm Cov}}

\newcommand{\fib}{\hbox{fib}}

\newcommand{\beg}{\mathop{{\rm beg}}}
\newcommand{\rank}{\mathop{{\rm rank}}}
\newcommand{\Loc}{\hbox{\rm Loc}}

\newcommand{\inv}{\mathop{{\rm inv}}}
\newcommand{\id}{\hbox{id}}

\newcommand{\D}{{D}}
\newcommand{\V}{{V}}

\newcommand{\Mon}{{\hbox{\rm  Mon}}}
\newcommand{\Aut}{{\hbox{\rm  Aut}}}

\makeatletter

\numberwithin{equation}{section}
\numberwithin{figure}{section}

\makeatother

\title[Two-orbit maniplexes]{An existence result on two-orbit maniplexes}
\author{Daniel Pellicer}
\author{Primo\v{z} Poto\v{c}nik}
\author{Micael Toledo}

\address{Micael Toledo,
University of Primorska, FAMNIT,
Glagolja\v{s}ka 8, Koper, Slovenia\newline
Also affiliated with: Institute of Mathematics and Physics, Jadranska 19, Ljubljana, Slovenia
} 
\email{micaelalexitoledo@gmail.com}

\address{Primo\v{z} Poto\v{c}nik, Faculty of Mathematics and Physics, University of Ljubljana, Jadranska 21, SI-1000 Ljubljana, Slovenia.\newline
\indent Also affiliated with: Institute of Mathematics, Physics and Mechanics, Jadranska 19, SI-1000 Ljubljana, Slovenia.
}
\email{primoz.potocnik@fmf.uni-lj.si}

\address{Daniel Pellicer\\
Research Centre for Mathematical Sciences, National Autonomous University of Mexico}
\email{pellicer@matmor.unam.mx}

\begin{document}
\begin{abstract}
A maniplex of rank $n$ is a connected, $n$-valent, edge-coloured graph that generalises abstract polytopes and maps. If the automorphism group of a maniplex $\mM$ partitions the vertex-set of $\mM$ into $k$ distinct orbits, we say that $\mM$ is a $k$-orbit $n$-maniplex. The symmetry type graph of $\mM$ is the quotient pregraph obtained by contracting every orbit into a single vertex. Symmetry type graphs of maniplexes satisfy a series of very specific properties. The question arises whether any pregraph of order $k$ satisfying these properties is the symmetry type graph of some $k$-orbit maniplex. We answer the question when $k = 2$.
\end{abstract}
\maketitle

\section{Introduction}
\label{sec:first}

A maniplex $\mM$ of rank $n$ can be defined as an $n$-regular properly edge-coloured graph satisfying some additional properties (see Section \ref{sec:intro}). It generalises a map on a surface to higher ranks, and at the same time is a relaxation of the notion of abstract polytope. Maniplexes were first introduced in \cite{SteveManiplex}, and later used for example in \cite{3Orb} and \cite{Twist}.

The symmetry of $\mM$ is measured by the number of orbits under the automorphism group,
where by an automorphism of $\mM$ we mean a color-preserving automorphism of the graph $\mM$.
 If this number is $k$ then $\mM$ is said to be a $k$-orbit maniplex. Under this notion, the most symmetric maniplexes are those where $k=1$, called {\em regular}.

The symmetry type graph of a maniplex $\mM$ is the quotient of $\mM$ by its automorphism group. For instance, the symmetry type graph of a regular maniplex of rank $n$ consists of one vertex and $n$ semi-edges, one of each colour. For a fixed $k \ge 2$, there may be several non-isomorphic candidates for symmetry type graphs of $k$-orbit maniplexes.

For example, when $k=2$ the possible symmetry type graphs of maniplexes of rank $n$ consist of two vertices, $\ell \ge 1$ edges between them, and $n-\ell$ semi-edges at each of the vertices. Since there are $2^n-1$ such properly $n$-coloured two-vertex graphs, there are at most $2^n-1$ symmetry type graphs of $2$-orbit maniplexes of rank $n$. 
Motivated by the terminology from the theories of maps and abstract polytopes, we call a $2$-orbit maniplex of rank $n$ whose symmetry type graph consists of two vertices and $n$ edges between them {\em chiral}. 

The existence of a chiral maniplex for each rank $n$ is a consequence of a result proved in \cite{Pell10}, which 
settled affirmatively a twenty-year-old question about existence of chiral abstract polytopes of every rank $n\ge 3$.

In this paper, we extend this consequence by showing that indeed each of the $2^n-1$ properly edge-coloured $n$-regular two-vertex graphs is the symmetry type graph of a $2$-orbit maniplex (see Theorem \ref{theo:main}).

As mentioned above, the maniplexes constructed in \cite{Pell10} also satisfy the necessary conditions to be abstract polytopes. This was proven using a vast background available on the automorphism groups of chiral polytopes. Nowadays there is still no analogous background for $2$-orbit polytopes. Therefore, determining whether the maniplexes constructed here are abstract polytopes or not would require much more work, and is left out of this paper.

In Sections \ref{sec:intro} and \ref{sec:bicolour} we recall basic concepts of maniplexes and polytopes. The polytope $\hat{2}^{\mM}$ is used in Section \ref{sec:2toM} to construct a family of polytopes that are the start point of our main construction. Then in Section \ref{sec:voltage} we relate maniplexes with the theory of graph coverings, and use this viewpoint to construct $2$-orbit $(n+1)$-maniplexes in Section \ref{sec:last}.

\section{Basic notions}
\label{sec:intro}

\subsection{Maniplexes}

A {\em maniplex} can be defined in several equivalent ways, one of them being as follows. For a positive integer $n$, an  {\em $n$-maniplex} $\mM$ is an ordered pair $(\mF(\mM),\{r_0,r_1,...,r_{n-1}\})$ where $\mF(\mM)$ is a non-empty set and $\{r_0,r_1,...,r_{n-1}\}$ is a set of fixed-point free involutory permutations of $\mF(\mM)$ satisfying the following three conditions: first, the permutation group $\langle r_0,r_1,...,r_{n-1} \rangle$ is transitive on $\mF$; second, for every $\Phi \in \mF$, its images under $r_i$ and $r_j$ are different if $i \neq j$; and third, for every $i,j \in \{0,...,n-1\}$ such that $|i-j| \geq 2$, the permutations $r_i$ and $r_j$ commute. The group $\langle r_0,r_1,...,r_{n-1} \rangle$  is then called the {\em monodromy group} of $\mM$ (also called the {\em connection group}) and is denoted $\Mon(\mM)$.
We shall call flags the elements of $\mF$, denote the image of a flag $\Phi$ under $r_i$ as $\Phi^i$ and say that $\Phi$ and $\Phi^{i}$ are $i$-adjacent. 

We can visualize an $n$-maniplex $(\mF,\{r_0,r_1,...,r_{n-1}\})$ as a simple connected graph with vertex-set $\mF$ and edge-set $\{\{\Phi,\Phi^i\}:\Phi \in \mF, i \in \{0,...,n-1\}\}$ together with a proper edge-colouring with the $i^{th}$ chromatic class consisting of edges $\{\{\Phi,\Phi^i\}:\Phi \in \mF\}$. Observe that the commutativity condition on the permutations $r_i$ implies that the induced subgraph with edges of colours $i$ and $j$, with $|i-j| \geq 2$, is a union of disjoint $4$-cycles.

Conversely, each connected $n$-regular graph with vertex-set $V$ and a proper edge colouring with chromatic classes $\{R_0,...,R_{n-1}\}$ satisfying the above condition on induced subgraphs arises from an $n$-maniplex $(V,\{r_0,r_1,...,n-1\})$, where $\Phi^i$ is the unique vertex such that $\{\Phi,\Phi^i\} \in R_i$. In this sense, an $n$-maniplex can be defined as a graph with an appropriate edge colouring. In what follows we will use both points of view interchangeably, as convenient.

Henceforth, let $\mM = (\mF,\{r_0,r_1,...,r_{n-1}\})$ be an $n$-maniplex.

For each $i \in \{0,1,...,n-1\}$ we define the set of $i$-faces of $\mM$ as the set of connected components of $\mM$ after the removal of all $i$-coloured edges. Equivalently, we may define an $i$-face as an orbit of the group $\langle r_j | j \in \{0, \ldots, n-1\} \setminus \{i\} \rangle$.

We may associate each $i$-face $F$ of $\mM$ with an $i$-maniplex by identifying two flags of $F$ whenever they are joined by a $j$-coloured edge, with $j \geq i$. Similarly as it is done with abstract polytopes, the $(n-1)$-faces of $\mM$ will be called the {\em facets} of $\mM$.

Let $p_i$ be the order of $r_{i-1} r_i$ for $i \in \{1, \dots, n-1\}$. Then the {\em Schl\"afli type} of $\mM$ is $\{p_1, \dots, p_{n-1}\}$. Whenever $Aut(\mM)$ acts transitively on the pairs $(F_{i-2},F_{i+1})$ consisting of an $(i-2)$-face $F_{i-2}$ and an $(i+1)$-face $F_{i+1}$ incident to $F_{i-2}$, $p_i$ indicates the number of $i$-faces incident simultaneously with $F_{i-2}$ and $F_{i+1}$. When $Aut(\mM)$ does not posses this transitivity property then $p_i$ is the least common multiple of these numbers of $i$-faces, when taking all possible pairs $(F_{i-2},F_{i+1})$. A regular map with $p$-gonal faces and degree $q$ has Schl\"afli type $\{p,q\}$.

An {\em automorphism} of the maniplex $\mM$ is a permutation of $\mF$ that
commutes with each $r_i$, $i\in \{0,\ldots,n-1\}$, and can thus be viewed as a colour-preserving automorphism of the associated coloured graph. Or equivalently, an automorphism of $\mM$ is a permutation of flags of $\mM$
 that maps each pair of $i$-adjacent flags to a pair of $i$-adjacent flags (this is stated for future reference in Lemma~\ref{lem:adj} below).
 As $\mM$ is connected, we see that the action on $\mF$ of its automorphism group, denoted $\Aut(\mM)$, is semiregular. Moreover, every automorphism is completely determined by its action on a single flag.

\begin{lemma}
\label{lem:adj}
Let $\mM$ be an $n$-maniplex and let $G$ be its automorphism group. Let $\Phi$ and $\Phi^i$ be two $i$-adjacent flags of $\mM$ and let $\Phi^G$ and $(\Phi^i)^G$ be their respective orbits under $G$. Then, for every flag $\Psi \in \Phi^G$ we have $\Psi^i \in (\Phi^i)^G$.
\end{lemma}

Recall that if 
$\Aut(\mM)$ has $k$ orbits of flags, we say that $\mM$ is a {\em $k$-orbit maniplex}. 
A $1$-orbit maniplex is called also called a {\em regular maniplex}, and if $\mM$ has two orbits of flags in such a way that adjacent flags belong to different orbits, the maniplex $\mM$ is said to be {\em chiral}.

Suppose $\Phi$ is a fixed flag of a regular $n$-maniplex $\mM$. Define $\rho_i$ as the unique automorphism of $\mM$ such that $\Phi^{\rho_i}=\Phi^i$. 
The following lemma is then a direct consequence of the commutativity of the action of $\Mon(\mM)$ and $\Aut(\mM)$.

\begin{lemma}
Let $\Phi$ be a fixed flag of $\mM$ and let $\alpha \in \Aut(\mM)$ be such that 
$\Phi^\alpha = \Phi^{r_{i_1} r_{i_2}\ldots r_{i_k}}$. Then, $\alpha= \rho_{i_k}\rho_{i_{k-1}}...\rho_{i_1}$. 
\end{lemma} 

It follows from the previous lemma and the fact that $\mM$ is connected that $\Aut(\mM)=\langle\rho_0,\rho_1,...,\rho_{n-1}\rangle$. We call $\{\rho_0,...,\rho_{n-1}\}$ 
the {\em standard generators} of $\Aut(\mM)$ relative to the base flag $\Phi$.

The dual of an $n$-maniplex $\mM = (\mF, \{r_0,r_1,\ldots,r_{n-1}\})$, denoted $\mM^*$, is the $n$-maniplex defined as $\mM^*=(\mF,\{r_0',r_1',\ldots,r_{n-1}'\})$, where $r_i'=r_n-i-1$. It is straightforward to see that the dual is well defined and that $(\mM^*)^*=\mM$. Moreover, $\Aut(\mM) \cong \Aut(\mM^*)$.\\

\subsection{Polytopes}

Here we briefly summarise some facts about abstract polytopes needed in this paper. We refer the reader to \cite{ARP} for
further information.
An {\em abstract $n$-polytope} is a partially ordered set $\mP$ with a monotone function $\rank:\mP \to \{-1,0,...,n\}$. Elements of $\mP$ of rank $i$ will be called {\em $i$-faces} and maximal chains in $\mP$ will be called {\em flags}. Two flags are said to be $i$-adjacent if they differ only in their $i$-face. We will also ask of $\mP$ to satisfy the following four conditions:

\begin{enumerate}
\item $\mP$ has a unique least face $F_{-1}$ and a unique greatest face $F_n$.
\item All flags in $\mP$ have length $n+2$.
\item 
For every $i \in \{0, \dots, n-1\}$ and every two faces $F$ and $G$, of ranks $i-1$ and $i+1$ respectively, such that $F<G$, there are exactly two $i$-faces $H_1$ and $H_2$ satisfying $F<H_1,H_2<G$. This is generally called the {\em diamond condition}. 
\item \label{c4}
If $\Phi$ and $\Psi$ are two flags, then there exists a sequence of successively adjacent flags $\Phi=\Phi_0,\Phi_1,...,\Phi_k=\Psi$ such that $\Phi \cap \Psi \subset \Phi_i$, for all $i \in \{0,...,k\}$. 

\end{enumerate}

We may associate each $i$-face $F$ of $\mP$ with the poset $\{G \in \mP: G \leq F\}$, which satisfy all the conditions to be an abstract polytope. In that sense, every $i$-face of $\mP$ is an abstract $i$-polytope in its own right. As is customary, $0$-faces will be called {\em vertices}, $1$-faces will be called {\em edges}, $2$-faces will simply be called faces
and $(n-1)$-faces will be called {\em facets}. 

Let $\mP$ be an abstract $n$-polytope and let $\mF$ be the set of all its flags. If $\Phi \in \mF$, then it follows easily from the diamond condition that there is a unique flag, denoted $\Phi^i$, that is $i$-adjacent to $\Phi$. For each $i \in \{0,...,n-1\}$, let $r_i$ be the permutation on $\mF$ that sends each flag to its $i$-adjacent flag, that is $\Phi^{r_i}=\Phi^i$. Define the {\em monodromy group} (also called  the {\em connection group})  of $\mP$ as $\Mon(\mP)=\langle r_0, r_1,...,r_{n-1}\rangle$ and note that, since $\mP$ satisfies the diamond condition, each generator $r_i$ is an involution. Furthermore, $r_ir_j=r_jr_i$ if $|i-j|\geq 2$. It follows from Condition (\ref{c4}) above  that the monodromy group of an abstract polytope acts transitively on its flags. This allows us to define an $n$-maniplex $\mM_\mP=(\mF, \{r_0,r_1,...,r_{n-1}\})$ {\em associated to $\mP$}. In this sense, an $n$-polytope can also be viewed as an $n$-maniplex.
The converse, however is not true: there are maniplexes that are not associated to any polytope (see, for example, \cite{SteveManiplex}).

An {\em automorphism} of the polytope $\mP$ is an order-preserving permutation of $\mP$, and the group of all automorphisms of $\mP$ is denoted by $\Aut(\mP)$. Since all the flags of $\mP$ have the same length, $\Aut(\mP)$ induces a faithful action on $\mF$. Moreover, every $\alpha \in \Aut(\mP)$ maps $i$-adjacent flags to $i$-adjacent flags, for all $i\in\{0,...,n-1\}$. It follows that the actions of $\Mon(\mP)$ and $\Aut(\mP)$ on $\mF$ commute. We say that an abstract polytope $\mP$ is regular if its automorphism group is transitive on flags and we say that $\mP$ is chiral if $\Aut(\mP)$ has two orbits of flags and adjacent flags always belong to different orbits. 
 
 Note that the notion of the monodromy group,  the automorphism group, the regularity and the chirality of a polytope coincide with those of the associated maniplex.

\subsection{Symmetry Type Graphs}

Intuitively speaking, the {\em symmetry type graph} of an $n$-maniplex $\mM$ is the
edge-coloured graph whose vertices are orbits of $\Aut(\mM)$ (in its action on the flags of $\mM$)
and with two orbits $A$ and $B$ adjacent by and edge of colour $i$ whenever there is a flag $\varphi\in A$
and a flag $\psi\in B$ that are $i$-adjacent in $\mM$. However, this definition requires the symmetry type graph to have edges that connect an orbit to itself. For that reason we define the following notion of a pregraph (which is widely used in the context of graph quotients and covers; see for example \cite{MNS}).

A {\em pregraph} is an ordered $4$-tuple $(\V,\D; \beg,\inv)$ where
$\V$ and $\D \neq \emptyset$ are disjoint finite sets of {\em vertices}
and {\em darts}, respectively, $\beg\colon \D \to \V$ is a mapping
which assigns to each dart $x$ its {\em initial vertex}
$\beg\,x$, and $\inv\colon \D \to \D$ is an involution which interchanges
every dart $x$ with its {\em inverse dart}, also denoted by $x^{-1}$.
The {\it neighbourhood} of a vertex $v$ is defined as the set
of darts that have $v$ for its initial vertex and the
{\it valence} of $v$ is the cardinality of the neighbourhood.

The orbits of $\inv$ are called {\em edges}.
The edge containing a dart $x$ is called a {\em semi-edge} if $x^{-1} = x$,
a {\em loop} if $x^{-1} \neq x$ while
$\beg\,(x^{-1}) = \beg\,x$,
and  is called  a {\em link} otherwise.
The {\em endvertices of an edge} are the initial vertices of the darts contained in the edge.
Two links are {\em parallel} if they have the same endvertices.

A pregraph with no semi-edges, no loops and no parallel links
can clearly be viewed as a (simple) graph. Conversely, every simple graph
 given in terms of its vertex-set $\V$ and edge-set $E$
can be viewed as the graph $(\V,\D; \beg,\inv)$, where
$D=\{(u,v) \mid uv \in E\}$,
$\inv (u,v) = (v,u)$ and  $\beg(u,v) = u$ for any $(u,v) \in D$. The mapping from $\Gamma$ to $\Gamma / G$ is called {\em the quotient projection with respect to $G$}.

An automorphism of a pregraph $\Gamma = (\V,\D; \beg,\inv)$ is a permutation
on $\D\cup \V$ that preserves both $\V$ and $\D$ and commutes with $\beg$ and $\inv$.
The group of all automorphisms of $\Gamma$ is denoted by $\Aut(\Gamma)$.

Given a pregraph $\Gamma = (\V,\D; \beg,\inv)$ and $G\le \Aut(\Gamma)$,
one can define the {\em quotient pregraph} $\Gamma/G = (\V/G, \D/G, \beg_G, \inv_G)$,
where $\V/G$ and $\D/G$ are the sets of all orbits of $G$ on $\V$ and $\D$, respectively, and
$\beg_G \colon \D/G \to \V/G$ and $\inv \colon \D/G \to \D/G$ are defined in such a way that
$\beg_G (x^G) = \beg(x)^G$ and $\inv_G (x^G) = \inv(x)^G$ for every $x\in \D$.

Let $\mM$ be a $k$-orbit $n$-maniplex viewed as an edge-coloured graph, in fact, as a coloured pregraph 
$\Gamma_\mM$ where each dart is coloured by the same colour as the underlying edge. Let $G$ be the group of colour preserving automorphism of $\Gamma_\mM$.
The {\em symmetry type graph} $T(\mM)$ of $\mM$ is then defined as the pregraph 
$\Gamma/G$, together with the colouring of the darts of $\Gamma/G$ where a dart $x^G$ is coloured by the same colour as the dart $x$ of $\Gamma_\mM$. To see that this colouring is well defined, consider a dart $x$ in $\Gamma_\mM$ and note that, since the elements of $G$ are colour-preserving, every  $y\in x^G$ must have the same colour as $x$, and so the colour of $x^G$ does not depend on its representative.
Furthermore, Lemma \ref{lem:adj} guarantees that, for each vertex $v^G$ of $T_{\mM}$, there is exactly one dart of each colour having $v^G$ as its initial vertex. It follows that, for each $i \in \{0,1,...,n-1\}$, there is exactly one edge
 of colour $i$ incident to each vertex of $T(\mM)$. In particular, this proves that
the symmetry type graph of a maniplex cannot have loops (and thus every edge is either a link or a semiedge).

Let $i,j \in \{0,...,n-1\}$ be such that $|i-j| \geq 2$. Since every walk of length $4$ in $\Gamma(\mM)$ of alternating colours $i$ and $j$ is necessarily a $4$-cycle, we have that any walk of length $4$ in $T(\mM)$ of alternating colours $i$ and $j$ must be a closed walk. We thus have that each connected component of the subgraph of $T(\mM)$ induced by edges (which could be links or semi-edges) of colour $i$ and $j$ is isomorphic to one of the following pregraphs. \\

\begin{figure}[h!]
\label{fig:stg}
\centering
\includegraphics[width=0.8\textwidth]{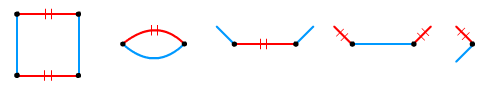}
\caption{Possible quotients of $2$-coloured $4$-cycles}
\label{fig:STGquo}
 \end{figure}

The symmetry type graph of a regular $n$-maniplex consist of a single vertex with $n$ semi-edges incident to it. Regular abstract polytopes have been widely studied (\cite{ARP} for details) and constructions of regular polytopes abound.

\begin{figure}[h!]
\centering
\includegraphics[width=0.8\textwidth]{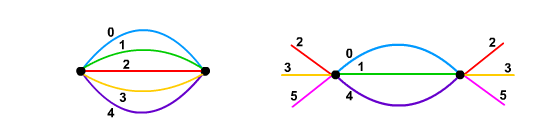}
\caption{The symmetry type graphs $2_{\emptyset}^5$ (left) and $2_{\{2,3,5\}}^6$ (right)}
\label{fig:STG}
 \end{figure}

The symmetry type graph $\Gamma$ of a $2$-orbit $n$-maniplex, is a connected $n$-valent, loopless pregraph of order $2$. Furthermore, the edges of $\Gamma$ are coloured in such a way that each vertex has either an edge or a semi-edge of each colour in $\{0,1,...,n-1\}$ incident to it. If $\mI$ is the set of colours appearing in the semi-edges of $\Gamma$, then $\Gamma$ is denoted by $2_{\mI}^n$ (see Figure \ref{fig:STG} for two examples). Following \cite{IsaThesis} and \cite{Isa10}, a maniplex having $2_{\mI}^n$ as its symmetry type graph is said to be {\em of type $2_{\mI}^n$}.

We are now ready to state our main theorem. Its proof is given in Section \ref{sec:last}.

\begin{theorem}
\label{theo:main}
For every $n \geq 3$ and every proper subset $\mI$ of $\{0,\ldots,n-1\}$ there exists an $n$-maniplex of type $2_{\mI}^n$.
\end{theorem}

Most of the works on polytopes and maniplexes has been devoted to regular ones \cite{ARP}. The most studied class of non-regular polytopes and maniplexes are the chiral ones. Its existence in ranks 3 and 4 was established in the 20-th century (see \cite{CoM} and \cite{twisted}), and in ranks 5 and higher only in the 21-th century (see \cite{higher} and \cite{Pell10}). Other classes of polytopes and maniplexes have been seldom studied \cite{korbit}. The study of the particular case of $2$-orbit polytopes and maniplexes started in \cite{IsaThesis} and \cite{Isa10}.

There were several attempts to show existence or non-existence of chiral polytopes of all ranks $n \ge 3$ in the years between 1991 (where the basic theory of chiral polytopes was established in \cite{chiral91}) and 2010 (when the existence was finally proven in \cite{Pell10}), and all faced great difficulties, particularly when ensuring that the constructed object is not regular. In view of this, it sounds a great task to determine if for every $n \ge 3$ and every $n$-regular pregraph $\Gamma$ satisfying the restrictions illustrated in Figure \ref{fig:STGquo} there is a maniplex with symmetry type $\Gamma$.

A more managable question is whether every graph $2_{\mI}^n$ is the symmetry type graph of some $2$-orbit maniplex, for each $n \geq 3$ and each $\mI \subsetneq \{0,1,...,n-1\}$. Since chiral polytopes exist for all ranks equal or greater than 3, every pregraph $2_{\emptyset}^n$ (one with no semi-edges) is the symmetry type graph of some (polytopal) maniplex.

We extend the work of the first author \cite{Pell10} by providing a construction that yields, for each $n \geq 3$ and each $\mI \subset \{0,1,...,n-1\}$, an $n$-maniplex of type $2_{\mI}^n$.

The solution of this problem only for pregraphs with two vertices may seem limited. However, the technique developed here may prove useful when attacking this problem for pregraphs with more vertices. We believe that the lack of symmetry on pregraphs (for instance, not being vertex-transitive) may be used to produce maniplexes that gain no extra symmetry.

Incidentally, it is known that every $n$-regular pregraph with $3$ vertices satisfying the restrictions illustrated in Figure \ref{fig:STGquo} is the symmetry type graph of a $3$-orbit $n$-polytope \cite{Cunn}. This is thanks to some hereditary property of the pregraphs with three vertices \cite{3Orb}. Then the present paper establishes that every $n$-regular pre-graph with at most $3$ vertices satisfying the restrictions illustrated in Figure \ref{fig:STGquo} is the symmetry type graph of a maniplex.

\section{Bi-colourings consistent with $\mI$}
\label{sec:bicolour}

Given an $n$-maniplex $\mM$ and $\mI \subseteq \{0,\dots, n-1\}$, a {\em bi-colouring of $\mM$ consistent with $\mI$} is a colouring of $\mF(\mM)$ with colours black and white such that $i$-adjacent flags have the same colour if and only if $i \in \mI$. In the context of maps by denoting $\mathcal{J}:=\{0,1,2\} \setminus \mI$, such a bi-colouring is called an $\mathcal{J}$-colouring in \cite{Wil2}.

For example, every $2$-maniplex is bipartite as a coloured graph, and hence it admits a bi-colouring consistent with $\emptyset$. In general, orientable maniplexes (that is, maniplexes that are bipartite as coloured graphs) admit a bi-colouring consistent with $\emptyset$. On the other hand, a triangle does not admit a bi-colouring consistent with $\{0\}$, since that would be equivalent to a proper $2$-colouring of the edges of the triangle.

The following proposition relates maniplexes that have bi-coulorings consistent to $\mI$ with cycles in the maniplex, in a similar way as bipartite graphs are related to the non-existence of odd cycles.

\begin{proposition}\label{pro:cycles}
Let $\mI \subseteq \{0, \dots, n-1\}$ and $\mM$ be an $n$-maniplex. Then $\mM$ admits a bi-colouring consistent with $\mI$ if and only if there is no cycle in $\mM$ (as a graph) having an odd number of edges whose labels are not in $\mI$.
\end{proposition}

\begin{proof}
A cycle in $\mM$ having an odd number of edges whose labels are not in $\mI$ cannot be bi-coloured in such a way that adjacent flags have the same colour only when $i \in \mI$. Hence, if $\mM$ contains such a cycle then it cannot admit a bi-couloring consistent with $\mI$.

Conversely, if all cycles in $\mM$ have an even number of edges whose labels are not in $\mI$, it is possible to assign a bi-couloring consistent with $\mM$ by couloring black any flag $\Phi$, and also any flag that can be reached from $\Phi$ by a path having an even number of edges whose labels are not in $\mI$. The remaining flags, namely those that can be reached from $\Phi$ by a path having an odd number of edges with labels that are not in $\mI$, are coloured white. This colouring is well-defined since otherwise it would be possible to construct a cycle in $\mM$ with an odd number of edges having labels not in $\mI$. Furthermore, this colouring is clearly consistent with $\mI$.
\end{proof}

When having a presentation for the automorphism group of a regular maniplex $\mM$, in terms of the standard generators, it is easy to verify if $\mM$ admits a bi-colouring consistent with some $\mI$. This is shown in the following proposition.

\begin{proposition}
\label{pro:bicolor}
Let $\mM$ be a regular maniplex of rank $n$, let $\Phi$ be a flag of $\mM$, let
$$
\Aut(\mM) = \langle \rho_0, \ldots, \rho_{n-1} \mid \mR \rangle
$$
be a presentation of $\Aut(\mM)$ in terms of the standard generators $\rho_0, \ldots, \rho_{n-1}$ (with respect to the base flag $\Phi$),
and let  $\mI \subseteq \{0,\ldots, n-1\}$. Suppose that for every $i\not \in \mI$ every relator in $\mR$ contains an even number of symbols $\rho_i$. Then $\mM$ admits a bi-colouring consistent with $\mI$.
\end{proposition}

\begin{proof}
Assume to the contrary that $\mM$ does not admit a bi-colouring consistent with $\mI$. Then, by Proposition \ref{pro:cycles} there is a closed cycle that has an odd number of edges whose labels are not in $\mI$. 
Since the standard generators $\rho_i$ were chosen with respect to the base flag $\Phi$, this can be rephrased as
\[\Phi=\Phi^{r_{i_1}r_{i_2}\dots r_{i_k}}=\Phi^{\rho_{i_k} \ldots \rho_{i_1}}\]
for some $r_{i_1},\dots,r_{i_k} \in \{r_0, \dots, r_{n-1}\}$, where $r_{i_j}\notin \mI$ for an odd number of $j$'s. It follows that $\rho_{i_k} \cdots \rho_{i_1}$ is a relator for $\Aut(\mM)$ with an odd number of symbols $\rho_i$ with $i \notin \mI$. Such a relator must be a product of conjugates of the relators in $\mR$, all of which have an even number of symbols $\rho_i$ with $i \notin \mI$. This contradics our hypothesis.
\end{proof}

\begin{remark}
\label{rem:allRelators}
Proposition \ref{pro:bicolor} requires a condition on one particular set $\mR$ of relators for $\Aut(\mM)$. However, if
this condition holds for one such set then it holds for every possible set of relators. Namely, if $w$ is a word consisting of symbols in $\{\rho_0, \ldots, \rho_{n-1}\}$ that equals $1$ in $\Aut(\mM) = \langle \rho_0, \ldots, \rho_{n-1} \mid \mR \rangle$ then it is contained in the normal closure $N$ in the group $\langle \rho_0, \ldots, \rho_{n-1} \rangle$ of the group generated by all conjugates of words in $\mR$. It is now obvious that if every relator contains an even number of symbols $\rho_i$ then so does every element in $N$.
\end{remark}

\section{Maniplex $\hat{2}^{\mM}$}\label{sec:2toM}

In this section we recall the construction of the structure $\hat{2}^{\mM}$. It was initially introduced by Danzer in \cite{Danzer}
in the context of incidence complexes, which are combinatorial structures more general than abstract polytopes.
This construction has appeared several times when studying regular polytopes and was generalised to maniplexes in \cite{Twist}.

Let $\mM = (\mF,\{r_0, \dots, r_{n-1}\})$ be an $n$-maniplex with set of facets $\mX$, $|\mX|\ge 2$.
 For $F\in \mX$ let $\chi_F \colon \mX \to \ZZ_2$ be the characteristic function of $F$ (that is, the function from $\mX$ to $\ZZ$ assigning $1$ to $F$ and $0$ to all other facets in $\mX$). Furthermore,
 for a flag $\Phi$ of $\mM$ let $F(\Phi)$ denote the facet of $\mM$ containing $\Phi$.

We now define the $(n+1)$-maniplex $\hat{2}^{\mM}=(\mF(\hat{2}^{\mM}),\{r_0', \dots, r_n'\})$ as follows:
\begin{eqnarray}\label{eq:Flags2M}
\mF(\hat{2}^{\mM}) &=& \mF(\mM) \times \ZZ_2^\mX,\\ \label{eq:Smallri}
(\Phi,x)^{r_i'} &=& (\Phi^{r_i},x) \mbox{ for } i < n,\\ \label{eq:Bigri}
(\Phi,x)^{r_n'} &=& (\Phi,x + \chi_{F(\Phi)}).
\end{eqnarray}

Note that $x + \chi_{F(\Phi)}$ in the equation above is the function in $\ZZ_2^\mX$ that differs from 
$x$ only when evaluated in the facet $F(\Phi)$.

\begin{figure}[h!]
\centering
\includegraphics[width=1.1\textwidth]{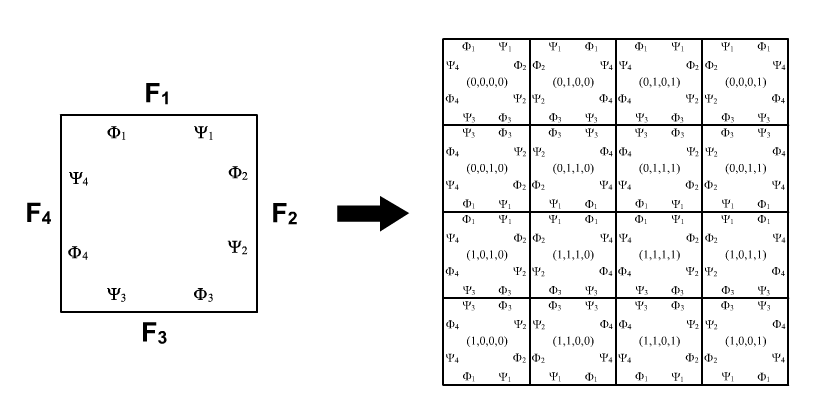}
\caption{The construction of $\hat{2}^{\mM}$ from a square $\mM$}
\label{fig:Delta}
 \end{figure}

The construction of the maniplex $\hat{2}^{\mM}$ from $\mM$ is illustrated in Figure \ref{fig:Delta}. On the left-hand side, we take$\mM$ to be the square with a set of facets $\mX=\{F_1,F_2,F_3,F_4\}$ and each facet $F_i$ containing flags $\Psi_i$ and $\Phi_i=\Psi_i^0$. The right-hand side depicts the maniplex $\hat{2}^{\mM}$ as a map on the torus, having $|\ZZ_2^\mX|$ copies of $\mM$ as facets, each labelled with an element of $|\ZZ_2^\mX|$, written as a vector. The glueing of the facets is given by (\ref{eq:Bigri}). For instance, facet $(0,1,1,0)$ is glued to $(0,1,1,1)$ through the edge with flags indexed $4$ because $(0,1,1,0) + \chi_{F_4} = (0,1,1,1)$.

For $\gamma\in\Aut(\mM)$ and $x \in \ZZ_2^\mX$, let the mapping
$x^\gamma: \mX \to \ZZ_2$ be defined by $x^\gamma=\gamma^{-1} \circ x$ and observe that this defines an embedding
of $\Aut(\mM)$ to the automorphism group of the elementary abelian $2$-group $\ZZ_2^\mX$.
This allows us to define the wreath product
\begin{equation}
\label{eq:wreath}
\ZZ_2 \wr \Aut(\mM) = \ZZ_2^\mX \rtimes \Aut(\mM).
\end{equation}

Furthermore, for $\gamma \in \Aut(\mM)$ and $y \in \ZZ_2^\mX$,
 let $\tilde{\gamma}$ and $\tilde{y}$ be the permutations
of $\mF(\hat{2}^\mM)$ defined by
\begin{eqnarray}
\label{eq:gammatilde}
(\Phi,x)^{\tilde{\gamma}} & = & (\Phi^\gamma, x^\gamma) \\
\label{eq:ytilde}
(\Phi, x)^{\tilde{y}} & = &(\Phi, x + y).
\end{eqnarray}
 
 Let $\Phi_0$ be the base flag of $\mM$ and $\rho_0, \dots, \rho_{n-1}$ the standard generators
 of $\Aut(\mM)$ with respect to $\Phi_0$.
Let $F_0$ be the base facet of $\mM$ (that is, the facet containing $\Phi_0$).
We choose the base flag of $\hat{2}^{\mM}$ to be $(\Phi_0, \overline{0})$, where
 $\overline{0}\in \ZZ_2^\mX$ is the $0$-constant function.
 
The following proposition summarises Section 6 of \cite{Twist}. Item (6) was proven in \cite{extensions}.

\begin{proposition}
\label{pro:oldresult}
Let $\mM$ be a regular maniplex with at least $2$ facets and Schl\"afli type $\{p_1, \dots, p_{n-1}\}$. Then the following hold:
\begin{enumerate}
 \item if $\gamma \in \Aut(\mM)$ then $\tilde{\gamma}\in \Aut(\hat{2}^\mM)$, 
 and if $y \in \ZZ_2^\mX$ then $\tilde{y} \in \Aut(\hat{2}^\mM)$;
 \item $\Aut(\hat{2}^{\mM})$ equals the wreath product $\ZZ_2 \wr \Aut(\mM)$ defined in (\ref{eq:wreath})
 and is generated by all $\tilde{\gamma}$ for $\gamma \in \Aut(\mM)$ and all $\tilde{y}$ for $y \in \ZZ_2^\mX$;
 \item  the stabiliser of the base facet of $\hat{2}^\mM$ equals 
 $\langle \tilde{\rho}_0, \ldots, \tilde{\rho}_{n-1}\rangle$ and there exists
 an isomorphism from $\Aut(\mM)$ to $\langle \tilde{\rho}_0, \ldots, \tilde{\rho}_{n-1}\rangle$
 mapping $\rho_i$ to $\tilde{\rho}_i$ for every $i\in \{0, \ldots, n-1\}$;
 \item
the standard generators of $\Aut(\hat{2}^{\mM})$ with respect to the base flag  $(\Phi_0, \overline{0})$
are $\tilde{\rho}_0, \dots, \tilde{\rho}_{n-1},\tilde{\chi_{F_0}}$ where $\chi_{F_0}$ is the characteristic function of the base facet $F_0$.
  \item $\hat{2}^{\mM}$ has Schl\"afli type $\{p_1, \dots, p_{n-1}, 4\}$;
 \item if $\mM$ is a polytope then so is $\hat{2}^{\mM}$.
\end{enumerate}
\end{proposition}

\begin{remark}\label{rem:not4}
All items of Proposition 6 hold if $\mM$ has only one facet, except for the fifth one. In this situation $\hat{2}^{\mM}$ coincides with the trivial maniplex over $\mM$ defined in \cite[Section 2.1]{Twist} and the last entry of the Schl\"afli type is $2$.
\end{remark}

\subsection{Symmetry type graph of $\hat{2}^{\mM}$}

Given a maniplex $\mM$, the next two lemmas describe transitivity properties of $\Aut(\hat{2}^{\mM})$. Throughout this subsection the Schl\"afli type of $\hat{2}^{\mM}$ plays no role, and Remark \ref{rem:not4} allows us to state these results for arbitrary maniplexes.

\begin{lemma}\label{le:FactTrans}
The automorphism group of the maniplex $\hat{2}^{\mM}$ acts transitively on the set of facets of $\hat{2}^{\mM}$.
\end{lemma}

\begin{proof}
Since a facet of $\hat{2}^{\mM}$ is an orbit of a flag under the group $\langle r_0, \dots, r_{n-1} \rangle$, from (\ref{eq:Smallri}) it follows that every facet of $\hat{2}^{\mM}$ can be described as
\[\{(\Phi,y_0) \, : \, \Phi \in \mF(\mM)\}\]
for some $y_0 \in \mathbb{Z}_2^{\mX}$. Then the elements $\tilde{y} \in \Aut(\hat{2}^{\mM})$ act transitively on the facets of $\hat{2}^{\mM}$ (see Proposition \ref{pro:oldresult} (1)).
\end{proof}

\begin{lemma}\label{le:2MkOrbit}
If $\mM$ is a $k$-orbit maniplex then so is $\hat{2}^{\mM}$.
\end{lemma}

\begin{proof}
Since all facets of $\hat{2}^{\mM}$ are isomorphic to $\mM$, the number of flag-orbits of $\hat{2}^{\mM}$ is at least $k$. On the other hand, the transitivity of $\Aut(\hat{2}^{\mM})$ on the set of facets and the automorphisms $\tilde{\gamma}$'s imply that $\hat{2}^{\mM}$ has precisely $k$ flag-orbits.
\end{proof}

The symmetry type graph of the maniplex $\hat{2}^{\mM}$ can be easily described in terms of that of $\mM$.

\begin{proposition}\label{pro:SymTypGra2}
Let $\mM$ be an $n$-maniplex. Then the symmetry type graph $T(\hat{2}^{\mM})$ is obtained from $T(\mM)$ by adding a semi-edge labelled $n$ at every vertex.
\end{proposition}

\begin{proof}
Any flag $\Phi$ is mapped to $\Phi^n$ by the automorphism $\tilde{\chi_{F(\Phi)}}$, implying that any two $n$-adjacent flags belong to the same flag-orbit. It follows that there is a semi-edge labelled $n$ at every vertex of $T(\hat{2}^{\mM})$.

By Lemma \ref{le:2MkOrbit}, $\mM$ and $\hat{2}^{\mM}$ have the same number of flag-orbits. Then, by Lemma \ref{le:FactTrans}, every flag-orbit has a representative in every facet. Since the facets of $\hat{2}^{\mM}$ are isomorphic to $\mM$, the flag-orbits of two flags $\Phi_1$ and $\Phi_2$ of $\mM$ are $i$-adjacent if and only if the flag orbits of the flags $(\Phi_1,y)$ and $(\Phi_2,y)$ of $\hat{2}^{\mM}$ are $i$-adjacent, for any $y \in \mathbb{Z}_2^{\mX}$. It follows that the graph obtained from $T(\hat{2}^{\mM})$ by removing all edges labelled $n$ is isomorphic to $T(\mM)$ as a labelled graph. This finishes the proof.
\end{proof}

\begin{corollary}\label{cor:edges0n}
Let $\mM$ be an $n$-maniplex with at least two facets and let $\mathcal{G}$ be an $(n+1)$-regular graph with edges labelled in $\{0, \dots, n\}$ that can be obtained from $T(\mM)$ by one of the two following procedures:
\begin{enumerate}
 \item Add semi-edges labelled $n$ at every vertex.
 \item Add $1$ to the label of every edge and add a semi-edge labelled $0$ at every vertex.
\end{enumerate}
Then $\mathcal{G}$ is the symmetry type graph of some maniplex $\mM'$. 
\end{corollary}

\begin{proof}
The symmetry type of the maniplexes $\hat{2}^{\mM}$ and $(\hat{2}^{\mM^*})^*$, where $\mM^*$ is the dual maniplex defined in 3.1, are $\mathcal{G}$ for the first and second procedures, respectively. 
\end{proof}

Similar arguments to the ones in this section can be used to prove Corollary \ref{cor:edges0n} using the trivial maniplex over $\mM$ defined in \cite{Twist} instead of $\hat{2}^{\mM}$.

\subsection{The polytopes $\mM_n$}

We now use the above construction to define recursively a family of polytopes.

\begin{definition}
\label{def:Mn}
Let $\mM_1$ be the unique regular polytope with rank $1$. Then, for $n \ge 2$ let $\mM_n$ be the regular polytope $\hat{2}^{\mM_{n-1}}$.
It is easy to see that $\mM_2$ is isomorphic to the square, and $\mM_3$ to the toroidal map $\{4,4\}_{(4,0)}$.

\end{definition}

The family of polytopes $\{\mM_n\}_{n \in \NN}$ has several interesting properties. One of them is given in the following proposition.

\begin{proposition}
For every $n \in \NN$ and $\mI \subseteq \{0,\dots,n-1\}$ the polytope $\mM_n$ admits a bi-colouring of flags consistent with $\mI$.
\end{proposition}

\begin{proof}
Let $G_n=\Aut(\mM_n)$ and let 
\begin{equation}
\label{eq:presentation}
G_n=\langle \rho_0, \ldots, \rho_{n-1} \mid \mR\rangle
\end{equation}
 be a presentation of $G_n$
in terms of standard generators. We will show that for every $i\in\{0,\ldots, n-1\}$ and every word $w\in \mR$
the symbol $\rho_i$ appears an even number of times in $w$.
Proposition~\ref{pro:bicolor} will then yield the result. By Remark~\ref{rem:allRelators} it suffices to
prove this condition for any set of relators $\mR$ for which (\ref{eq:presentation}) holds.

The proof of the above fact is by induction on $n$.
For $n=1$ we may choose a presentation of $G_1$ to be
$G_1= \langle \rho_0 \mid \rho_0^2\rangle$, which clearly satisfies the condition above.

Now assume that the claim holds for $G_{n-1}$. Let
\begin{equation}
\label{eq:presentation-1}
G_{n-1}=\langle \rho_0, \ldots, \rho_{n-2} \mid \mR' \rangle
\end{equation}
be a presentation for $G_{n-1}$ in terms of standard generators.
Let $\mX$ be the set of facets of $\mM_{n-1}$ and fix a base facet $F_0\in \mX$.
By Proposition~\ref{pro:oldresult} (2), $G_n = \ZZ_2^\mX \rtimes G_{n-1}$ where $G_{n-1}$ acts upon 
the elementary abelian group $\ZZ_2^\mX$ according to the rule $x^\gamma = \gamma^{-1} \circ x$ for every
$x\in \ZZ_2^\mX$ and $\gamma\in G_{n-1}$. Furthermore, let $e_0 = \chi_{F_0} \in \ZZ_2^\mX$
 be the characteristic function of the base facet $F_0$.

Recall that the standard generators of $G_n$ are $\tilde{\rho}_0, \ldots, \tilde{\rho}_{n-2}, \tilde{e}_0$ and
that the mapping $\gamma\mapsto \tilde{\gamma}$ in an embedding of $G_{n-1}$ to $G_n$.

Let $w$ be a word in the symbols $\{\tilde{\rho}_0, \ldots, \tilde{\rho}_{n-2}, \tilde{e}_0\}$ that evaluates
to $1$ in $G_n$. Then $w$ can be written as
$$w
= \tilde{\gamma}_1 x_1 \tilde{\gamma}_2 x_2 \cdots x_{k-1} \tilde{\gamma}_k
$$ where $\tilde{\gamma}_i$'s are words in $\{\tilde{\rho}_0, \ldots, \tilde{\rho}_{n-2}\}$
and $x_i = \tilde{e}_0^{s_i}$ for some integer $s_i$.
Observe that in $G_n$ the word $w$ evaluates to
$$
\tilde{\gamma}_1 \tilde{\gamma}_2 \cdots \tilde{\gamma}_k \, 
x_1^{\tilde{\gamma}_2  \cdots \tilde{\gamma}_k} 
x_2^{\tilde{\gamma}_3  \cdots \tilde{\gamma}_k} 
\cdots x_{k-1}^{\tilde{\gamma_k}}.
$$
Hence, in $G_n$ the following holds:
$$
\tilde{\gamma}_1 \tilde{\gamma}_2 \cdots \tilde{\gamma}_k \, 
x_1^{\tilde{\gamma}_2  \cdots \tilde{\gamma}_k} 
x_2^{\tilde{\gamma}_3  \cdots \tilde{\gamma}_k} 
\cdots x_{k-1}^{\tilde{\gamma_k}}  = 1, 
$$
which implies that
$$
\tilde{\gamma}_1 \tilde{\gamma}_2 \cdots \tilde{\gamma}_k = 1 \>\> \hbox{ and } \>\>
x_1^{\tilde{\gamma}_2 \cdots \tilde{\gamma}_k} \cdots x_{k-1}^{\tilde{\gamma_k}}  = 1. 
$$
But then
\begin{equation}
\label{eq:gamma}
\gamma_1 \gamma_2 \cdots \gamma_k
\end{equation}
is a trivial element of $G_{n-1}$,  and by induction hypothesis, for every $i\in \{0,\ldots,n-2\}$ each
symbol $\rho_i$ appears an even number of times in (\ref{eq:gamma}). Consequently, each 
symbol $\tilde{\rho}_i$ appears an even number of times in $\tilde{\gamma}_1 \tilde{\gamma}_2 \cdots \tilde{\gamma}_k$. Moreover, since each $x_j$ is a power of $\tilde{e}_0$, the symbols
$\tilde{\rho}_i$ also appear an even number of times in the expression
$
x_1^{\tilde{\gamma}_2 \cdots \tilde{\gamma}_k} \cdots x_{k-1}^{\tilde{\gamma_k}}.
$
Now let $S$ be the set of indices $i$ for which $s_i$ is odd. Since $\tilde{e}_0^2 = 1$
and since $x_1^{\tilde{\gamma}_2 \cdots \tilde{\gamma}_k} \cdots x_{k-1}^{\tilde{\gamma_k}}  = 1$, we see that
$$
\prod_{i\in S} \tilde{e}_0^{\tilde{\gamma}_{i+1} \cdots \tilde{\gamma}_k}  = 1.
$$
However, a factor $\tilde{e}_0^{\tilde{\gamma}_{i+1} \cdots \tilde{\gamma}_k}$, $i\in S$, is a characteristic
function of the facet $F_0^{\tilde{\gamma}_{i+1} \cdots \tilde{\gamma}_k}$, implying that the size of $S$
is even. This shows that $\tilde{e}_0$ appears an even number of times in the expression
$x_1^{\tilde{\gamma}_2 \cdots \tilde{\gamma}_k} \cdots x_{k-1}^{\tilde{\gamma_k}}$.
\end{proof}

The following lemmas show that $\mM_m$ satisfies some technical conditions required in Section \ref{sec:last} for the main construction.

\begin{lemma}
Let $\mM$ be a regular $n$-maniplex, and let $\mX_0$ be a set of facets of $\mM$ such that $\mX_0$ is not invariant under any non-trivial element of $\Aut(\mM)$. Then there exists an involutory element $\eta \in \Mon(\hat{2}^\mM)$ such that it maps every two distinct flags of $\hat{2}^{\mM}$ in the same facet to flags in distinct facets.
\end{lemma}

\begin{proof}
We choose the base facet $F_0$ of $\mM$ in such a way that $F_0 \in \mX_0$. For every $F \in \mX_0$ we choose a flag $\Phi_F$, and we let $\Phi_0 = \Phi_{F_0}$ be the base flag. For each flag $\Phi_F$ let $\omega_F$ be the element in $\Mon(\mM)$ that maps $\Phi_0$ to $\Phi_F$, and let $\mW = \{\omega_F \colon F \in \mX_0\}$.

Denote the standard generators of $\mM$ by $r_0, \dots, r_{n-1}$ and let the standard generators $r_0', \dots, r_n'$ of $\hat{2}^{\mM}$ be as in (\ref{eq:Smallri}) and (\ref{eq:Bigri}). Note that for every $\omega \in \mW$, the conjugate $(r_n')^\omega$ maps a flag $(\Phi,x)$ of $\hat{2}^\mM$ to $(\Phi, x + \chi_{F(\Phi^\omega)})$.
Let $\eta=\prod_{\omega \in \mW} (r_n')^\omega \in \Mon(\hat{2}^{\mM})$, that is,
$$
\eta : (\Phi, x) \mapsto (\Phi, x + \sum_{\omega \in \mW} \chi_{F(\Phi^\omega)}).
$$
Clearly, $\eta$ is an involution.

Let $\Phi$, $\Psi$ be two distinct flags in the same facet of $\hat{2}^{\mM}$. By regularity of $\hat{2}^{\mM}$ we may assume that $\Phi = \Phi_0$, and $\Psi$ can be expressed as $\Phi^\nu$ for some non-trivial $\nu \in \langle r_0', \dots, r_{n-1}' \rangle$. It suffices to show that $\Phi_0^\eta$ and $\Phi_0^{\nu \eta}$ belong to distinct facets, that is, there exists no $\nu' \in \langle r_0', \dots, r_{n-1}' \rangle$ such that $\Phi_0^{\nu \eta} = \Phi_0^{\eta \nu'}$. Since the action of the monodromy group of $\hat{2}^{\mM}$ is regular on the flags, this is equivalent to saying that $\nu^\eta \notin \langle r_0', \dots, r'_{n-1} \rangle$.

We know that
\begin{eqnarray}  (\phi,x)^{\nu^{\eta}}&=&(\phi,x+\sum_{\omega\in\mW}\chi_{F(\phi^\omega)})^{\nu\eta} = (\phi^\nu,x+\sum_{\omega\in\mW}\chi_{F(\phi^\omega)})^{\eta}\\ \label{eq:eqsets} &=& (\phi^\nu,x+\sum_{\omega\in\mW}(\chi_{F(\phi^\omega)}+\chi_{F(\phi^{\nu\omega})})).  
\end{eqnarray} 

From the definition of the standard generators of $\Mon(\hat{2}^\mM)$, the element $\nu^\eta$ belongs to $\langle r_0', r_1', \dots, r_{n-1}'\rangle$ if and only if $\sum_{\omega\in\mW}(\chi_F(\phi^\omega)+\chi_F(\phi^{\nu\omega}))$ is trivial.

Write $\nu=r_{i_1}'r_{i_2}' \cdots r_{i_k}'$ and let $\{\rho_0,\rho_1,\dots,\rho_{n-1}\}$ be the standard generators of $\Aut(\mM)$. Then $$\phi^{\nu\omega}=\phi^{r_{i_1}'r_{i_2}' \cdots r_{i_k}'\omega}=(\phi\rho_{i_k} \cdots \rho_{i_1})^\omega=\phi^\omega\rho_{i_k} \cdots \rho_{i_1}$$
and therefore $F(\phi^{\nu\omega})=F(\phi^\omega\rho_{i_k} \cdots \rho_{i_1})=F(\phi^\omega)\rho_{i_k} \cdots \rho_{i_1}$.

For every $\omega\in\mW$ the facet $F(\phi^\omega)$ belongs to $\mX_0$, then the expression in \ref{eq:eqsets} can be written as 
$$ (\phi^\nu,x+\sum_{F\in\mX_0}\chi_F + \sum_{F\in\mX_0 \rho_{i_k} \cdots \rho_{i_1}}  \chi_F).$$ 
By hypothesis, $\mX_0 \rho_{i_k} \cdots \rho_{i_1} \neq \mX_0$, and therefore $\sum_{\omega\in\mW}(\chi_F(\phi^\omega)+\chi_F(\phi^{\nu\omega}))$ is non-trivial and $\nu^\eta$ does not belong in $\langle r_0', r_1', \dots, r_{n-1}'\rangle$.
\end{proof}

\begin{figure}[h!]
\centering
\includegraphics[width=0.4\textwidth]{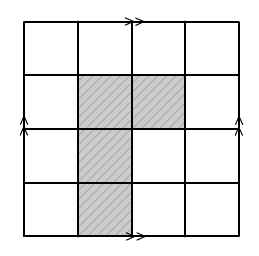}
\caption{The toroidal map $\{4,4\}_{(4,0)}$. The set of shaded facets is fixed only by the identity}
\label{fig:torusfacets}
 \end{figure}

\begin{lemma}
Let $m$ be an integer, $m\ge 3$, and let $\mM=\mM_m$ be the regular $m$-maniplex defined in Definition~\ref{def:Mn}.
Then there exists a set of facets $\mX_0$ of $\mM$ such that $\mX_0$ is not invariant under any non-trivial element of $\Aut(\mM)$. 
\end{lemma}

\begin{proof}
The proof is by induction on $m$. If $m=3$, then $\mM$ is isomorphic to the toroidal map  $\{4,4\}_{(4,0)}$ shown in Figure~\ref{fig:torusfacets},
 viewed as a $3$-maniplex.
The set $\mX_0$ can then be chosen to be the set of grey faces in the figure. The reader can convince herself that $\mX_0$ is fixed by no non-trivial automorphism of the map $\mM$.

Now let $n$ be an integer, $n\ge 3$, and suppose that the lemma holds with $m=n$.
We will show that it then also holds for $m=n+1$.
Let $\mM = \mM_{n}$, so that $\mM_m = \mM_{n+1}=\hat{2}^{\mM}$.
Further, let $\mX$ be the set of facets of $\mM$ and let
$\mX_0 \subset \mX$ be a set of facets fixed by no non-trivial automorphism of 
$\mM$.

Recall that a facet $F$ of an $(n+1)$-maniplex determined by the monodromy group
$\langle r_0, \ldots, r_{n}\rangle$ equals an orbit of $\langle r_0, \ldots, r_{n-1}\rangle$
acting on the set of flags of the maniplex. 

Further, recall that the set of flags of 
$\hat{2}^{\mM}$ is $\mF(\mM) \times \ZZ_2^\mX$ (see (\ref{eq:Flags2M})) and that
the monodromy group of $\hat{2}^{\mM}$ equals $\langle r_0', \ldots, r_n'\rangle$
where
$(\Phi,x)^{r_i'} = (\Phi^{r_i},x) \mbox{ for } i < n$ (see (\ref{eq:Smallri})) and $r_n'$ is as in (\ref{eq:Bigri}).

 In particular, a typical facet of
$\hat{2}^\mM$ equals the set
$$
 F_x = \{ (\Phi,x) \colon \Phi \in \mF(\mM) \}
$$
where $x$ is an arbitrary element of $\ZZ_2^\mX$.

We will now define a set of facets $\mX_0'$ of $\hat{2}^\mM$ and prove that it
is fixed by no non-trivial automorphism of $\hat{2}^\mM$:
$$
\mX_0' = \{F_{\chi_H} : H \in S_0\} \cup \{F_{\bar{0}}\}.
$$
(recall that $\bar{0}$ denotes the constant $0$ function).

We claim that the facet $F_{\bar{0}}$ is the only facet in $\mX_0'$ such that
 for each $H\in S_0$ it
contains a flag which is $n$-adjacent to some
flag in $F_{\chi_H}$. This will then imply that every automorphism of $\hat{2}^\mM$
that preserves $\mX_0'$ fixes the facet $F_{\bar{0}}$.

To prove the claim, note that if $\Phi$ is a flag in $H$ then $(\Phi,\bar{0})^{r_n}$ belongs to $F_{\chi_H}$,
while $(\Phi,\bar{0})$ belongs to $F_{\bar{0}}$. Now let $H, H' \in \mX$.
Then a flag in $F_{\chi_H}$ is of the form $(\Phi,\chi_H)$, while a flag in $F_{\chi_{H'}}$ is of the form $(\Phi',\chi_{H'})$ for some $\Phi \in H$ and $\Phi' \in H'$. Since $\chi_H$ and $\chi_{H'}$ differ in precisely two values, these two flags cannot be $n$-adjacent in $\bar{2}^\mM$.
This completes the proof of our claim that every automorphism of $\hat{2}^\mM$
that preserves $\mX_0'$ fixes the facet $F_{\bar{0}}$.

By part (2) of Proposition~\ref{pro:oldresult} an automorphism of $\hat{2}^\mM$
can be written uniquely as a product $\tilde{y}\tilde{\gamma}$ for some
$\gamma\in\Aut(\mM)$ and $y\in\ZZ_2^\mX$, as defined in (\ref{eq:gammatilde}) and
(\ref{eq:ytilde}); that is
$$
(\Phi,x)^{\tilde{\gamma}}  =  (\Phi^\gamma, x^\gamma) \> \hbox{ and } \>
(\Phi, x)^{\tilde{y}}  = (\Phi, x + y).
$$

 Suppose that $\tilde{y}\tilde{\gamma}$ preserves $\mX_0'$. We have shown that
then $\tilde{y}\tilde{\gamma}$ fixes $F_{\bar{0}}$. Note that $\bar{0}^\gamma = \bar{0}$
and hence
$$
F_{\bar{0}}^{\tilde{\gamma}} = \{ (\Phi^\gamma,\bar{0}^\gamma) \colon \Phi \in \mF(\mM) \}
= F_{\bar{0}}.
$$
This implies that $\tilde{y}$ fixes $F_{\bar{0}}$, which clearly implies that $y=\bar{0}$
and thus $\tilde{y}$ is trivial. 
Further, observe that $(F_{\chi_H})^{\tilde{\gamma}} = F_{\chi_{(H^{\gamma})}}$,

\begin{eqnarray*}
(F_{\chi_H})^{\tilde{\gamma}} &=& \{(\Phi,\chi_H) \, : \, \Phi \in \mF(\mM)\}^{\tilde{\gamma}}\\
&=& \{(\Phi^{\gamma},(\chi_H)^\gamma) \, : \, \Phi \in \mF(\mM)\} \\
&=& \{(\Phi,\gamma^{-1} \circ \chi_H) \, : \, \Phi \in \mF(\mM)\}= F_{\chi_{H^\gamma}},
\end{eqnarray*}
implying that
$$
\mX_0' = \mX_0'^{\tilde{\gamma}} = \{(F_{\chi_H})^{\tilde{\gamma}} \colon H\in \mX_0\} 
\cup \{F_{\bar{0}}^{\tilde{\gamma}}\} =
\{F_{\chi_{(H^{\gamma})}} \colon H\in \mX_0\} \cup \{F_{\bar{0}}\}.
$$
By the inductive hypothesis $\id$ is the only automorphism of $\mM$
that preserves $\mX_0$, hence $\gamma = \id$ and the lemma follows.

\end{proof}

\section{Voltage construction}\label{sec:voltage}

In this section we briefly summarise a few facts about covering projections and voltage assignments that are
significant for the proof of our main theorem. We refer the reader to \cite{MNS} for further background on this topic.
We then relate this theory with the notion of a $2$-orbit maniplex.

Let $\Gamma= (\D,\V; \beg,\inv)$ and $\Gamma'= (\D',\V'; \beg',\inv')$ be two pregraphs.
A {\em morphism of pregraphs}, $f \colon \Gamma \to \Gamma'$,
is a function $f \colon \V \cup \D \to \V' \cup \D'$
such that
$f(\V) \subseteq \V'$, $f(\D) \subseteq \D'$,
$f\circ \beg = \beg' \circ f$ and $f \circ \inv = \inv' \circ f$.
A graph morphism is an {\em epimorphism} ({\em automorphism}) if it is
a surjection (bijection, respectively).
Furthermore, $f$ is a {\it covering projection} provided that it is an epimorphism
and if for every $v\in \V$ the restriction of $f$ onto the neighbourhood of $v$ is a bijection.
If $f$ is a covering projection, then the group of automorphisms, that for every $x\in \V'\cup \D'$
preserve the preimage $f^{-1}(x)$ setwise, is called the group of covering transformations.
If, in addition, the graph $\Gamma$ is connected and the group of covering transformations is
transitive on each such preimage, then the covering projection $f$ is said to be {\it regular}
and the graph $\Gamma$ is then called a {\it regular cover} of $\Gamma'$.

Note that if $H \le Aut(\Gamma)$ acts semiregularly on $\V(\Gamma)$, then the resulting quotient projection from $\Gamma$ to $\Gamma / H$ is a covering projection.

Let $\Gamma' = (\D',\V',\beg', \inv')$ be an arbitrary connected pregraph, let 
$N$ be a group and let $\zeta \colon \D' \to N$ be a  mapping (called a 
{\em voltage assignment}) satisfying the condition 
$\zeta(x) = \zeta(\inv' x)^{-1}$ for every $x\in \D'$. 
Then $\Cov(\Gamma',\zeta)$ is the graph with $\D' \times N$ and $\V' \times N$ 
as the sets of darts and vertices, respectively, and the functions 
$\beg$ and $\inv$ defined by $\beg(x,a) = (\beg' x,a)$ and 
$\inv(x,a) = (\inv' x, \zeta(x) a )$. If, for a voltage assignment $\zeta\colon \D(\Gamma') \to N$,
there exists a a spanning tree $T$ in $\Gamma'$ such that a $\zeta(x)$ is the trivial element of $N$
for every dart in $T$, then we say that $\zeta$ is {\it normalised};   
note that in this case $\Cov(\Gamma',\zeta)$ is connected if and only if the image of $\zeta$ generates $N$.

It is well known that every regular cover of a given (pre)graph $\Gamma'$ is isomorphic to $\Cov(\Gamma',\zeta)$ for some
normalised voltage assignment $\zeta$. To be precise,
let $\wp \colon \Gamma \to \Gamma'$ be a regular covering projection and let $T$ be a spanning tree in $\Gamma'$. Then
there exists a voltage assignment $\zeta\colon \D(\Gamma') \to N$ which is normalised with respect to $T$,
such that $\Gamma \cong \Cov(\Gamma',\zeta)$.

We define a walk of length $n$ in a pregraph as a sequence $(x_1,x_2,..,x_n)$ of darts such that $\beg(x_{i+1})=\beg(x_i^{-1})$ for all $i \in \{1,...,n-1\}$. A walk is said to be {\it reduced} if it contains no subsequence of the form $x_i,x_i^{-1}$. From this point on all walks will be reduced and will be referred to simply as walks. 

If $\zeta \colon \D(\Gamma) \to G$ is a voltage assignment, then
we define the {\em net voltage} of a walk $(x_1,x_2,...,x_n)$ in $\Gamma$ as the product $\zeta(x_n)\zeta(x_{n-1})\ldots\zeta(x_1)$
in the voltage group. We will slightly abuse notation and use $\zeta(W)$ to denote the net voltage of a walk $W$.
For a vertex $u$, the {\em local voltage group} $\Loc(G,u)$ is the group of net voltages of all closed walks based at $u$. The fiber of $u$ under a covering projection $\wp: \Cov(\Gamma, \zeta)$ will be denoted by $\fib_u$.

Let us now move to the special situation of symmetry type graphs $\Gamma$ on two points introduced in Section~\ref{sec:intro}. Recall that the darts of such a pregraph $\Gamma$
 are coloured in such a way that each vertex has either a link or a semi-edge
  of each colour in $\{0,1,...,n-1\}$ incident to it (see figure \ref{fig:stg}).

Let $\Gamma$ be the pregraph $2_{\mI}^n$ and let $u$ and $v$ be its vertices. Denote by $u_i$ the unique $i$-coloured dart having initial vertex $u$. Define $v_i$ similarly. Observe that $\inv(u_i)=u_i$ whenever $i \in \mI$ and  $\inv(u_i)=v_i$ otherwise. The following two propositions give sufficient conditions on a voltage assignment $\zeta$ for $\Gamma$ to lift into an $n$-maniplex.

\begin{proposition}
\label{prop:manilift1}
Let $n\ge 3$ be an integer, let $\mI \subseteq \{2,\ldots, n-1\}$, and let $\Gamma$ be the pregraph $2_{\mI}^n$.
 Let $G$ be a group with distinct non-trivial generators $\{z,y_2,y_3,...,y_{n-1}\}$ such that
 \begin{eqnarray}\label{eq:volt1}
y_i^2 &=& 1_G \\ \label{eq:volt2}
z^{y_i} &=& \left\{
\begin{array}{ll}
z & \mbox { if } i\ge 3 \mbox{ and }  i \in \mI \\
z^{-1} & \mbox{ if } i\ge 3 \mbox{ and } i \notin \mI
\end{array}\right.\\ \label{eq:volt3}
[y_i,y_j] &=& 1_G \quad \mbox{ if }\> |i-j|>1 
\end{eqnarray}
If $\zeta : D(\Gamma) \to G$  is the voltage assignment satisfying $\zeta(u_0)= 1_G$, $\zeta(u_1)=z$ and $\zeta(u_i)= \zeta(v_i)=y_i$ for $i \geq 2$, then $\Cov(\Gamma,\zeta)$  is a maniplex.
\end{proposition}

\begin{proof}
Let $\mM:=\Cov(\Gamma,\zeta)$. Now, in $\mM$, colour every dart in the fiber of $x$ with the same colour $x$ has in $\Gamma$. Note that $\zeta$ is normalised with respect to the spanning tree induced by the dart $u_0$, and therefore $\mM$ is connected. Moreover, it is straightforward to see that $\mM$ is an $n$-valent simple graph. It remains to show that the subgraph induced by colour $i$ and $j$ is a union of squares whenever $|i-j| \geq 2$.  
Clearly, a closed walk in $\Gamma$ lifts to a closed walk in $\mM$ if and only if its net voltage is $1_G$. By observing $\Gamma$, it is not too difficult to notice that any walk of length $4$ in $\mM$ that traces darts of only two colours $i$ and $j$ with $|i-j|\ge 2$ necessarily ends and starts in vertices belonging to the same fibre. By taking a closer look at $\Gamma$, we also see that any such walk $W$ has net voltage $1_G$; for example: 
if $i \geq 2$ and $j = 0$, the net voltage of $W$ is either $y_i (1_G) y_i (1_G)$ or $y_i^{-1} (1_G) y_i^{-1} (1_G)$. From \ref{eq:volt1} we get that this equals $1_G$ in both cases.
Further, suppose that $i \geq 3$ and $j = 1$. If $i \in \mI$, then  the net voltage $W$ is either $y_i z y_i z^{-1}$ or $y_i z^{-1} y_i z$. From \ref{eq:volt2} and the fact that $y_i$ is an involution, we have that $y_i z y_i z^{-1} = (y_i z y_i) z^{-1} = zz^{-1}= 1_G$ and $y_i z^{-1} y_i z= (y_i z^{-1} y_i) z = z^{-1}z = 1_G$. Now, if $i \notin \mI$, then $W$ has net voltage $y_i z^{-1} y_i z^{-1}$  or  $y_i z y_i z$. Again, using \ref{eq:volt1} and \ref{eq:volt2}, we get $y_i z^{-1} y_i z^{-1} = (y_i z^{-1} y_i) z^{-1}=zz^{-1}= 1_G$ and $y_i z y_i z = z^{-1}z= 1_G$.
Finally, if $i,j \geq 2$, the result follows from \ref{eq:volt1} and \ref{eq:volt3}.
Since every closed walk of length $4$ in a simple graph is a cycle, this implies that every walk of length $4$ in $\mM$ with alternating colour $i$ and $j$, where  $|i-j| \geq 2$, is a $4$-cycle.
\end{proof}

\begin{proposition}
\label{prop:manilift2}
Let $n\ge 3$ be an integer, let $\mI \subseteq \{1,\ldots, n-1\}$ with $1 \in \mI$, and let $\Gamma$ be the pregraph $2_{\mI}^n$.
 Let $G$ be a group with generators $\{z,z',y_2,y_3,...,y_{n-1}\}$ such that
\begin{eqnarray}\label{eq:1inivolt1}
z^2 &=& (z')^2 = y_i^2 = 1_G \\ \label{eq:1inivolt2}
z^{y_i} &=& \left\{
\begin{array}{ll}
z & \mbox { if } i\ge 3 \mbox{ and }  i \in \mI \\
z' & \mbox{ if } i\ge 3 \mbox{ and } i \notin \mI
\end{array}\right.\\ \label{eq:1inivolt3}
(z')^{y_i} &=& \left\{
\begin{array}{ll}
z' & \mbox { if } i\ge 3 \mbox{ and }  i \in \mI \\
z & \mbox{ if } i\ge 3 \mbox{ and } i \notin \mI
\end{array}\right.\\ \label{eq:1inivolt4}
[y_i,y_j] &=& 1_G \quad \mbox{ if } |i-j|>1 
\end{eqnarray}
If $\zeta : D(\Gamma) \to G$  is the voltage assignment satisfying $\zeta(u_0)=1_G$, $\zeta(u_1)=z$, $\zeta(v_1)=z'$ and $\zeta(u_i)= \zeta(v_i)= y_i$ for $i \geq 2$, then $\Cov(\Gamma,\zeta)$  is a maniplex.
\end{proposition}

\begin{figure}[h!]
\centering
\includegraphics[width=1\textwidth]{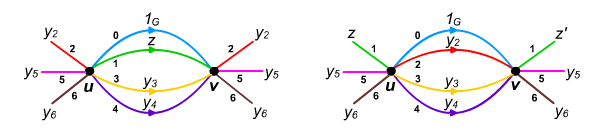}
\caption{Voltage assignment for $2_{\{2,5,6\}}^7$ and $2_{\{1,5,6\}}^7$.}
\label{fig:STGvolt}
 \end{figure}

In Figure \ref{fig:STGvolt} we show examples of the voltage assignments described in Propositions \ref{prop:manilift1} and \ref{prop:manilift2}.

\begin{proof}
The proof is almost identical to the proof of Proposition \ref{prop:manilift1}, except for the case when $i \geq 3$ and $j = 1$ which, as the reader can verify, follows from \ref{eq:1inivolt2}, \ref{eq:1inivolt3} and the fact that both $z$ and $z'$ are involutions.
\end{proof}

Suppose now that $G$ satisfies the relations \ref{eq:volt1}, \ref{eq:volt2} and \ref{eq:volt3} in Proposition~\ref{prop:manilift1} (if $1\not \in \mI$) or relations \ref{eq:1inivolt1}, \ref{eq:1inivolt2}, \ref{eq:1inivolt3} and \ref{eq:1inivolt4} in Proposition~\ref{prop:manilift2} (if $1\in \mI$), so that $\Cov(G,\zeta)$ is a maniplex.
Observe that $G$ acts on the vertices of $\Cov(\Gamma,\zeta)$ by right multiplication on the second coordinate. This is, for every $g\in G$, we can define a mapping $f_g$ given by $(u,h) \mapsto (u,hg)$. It is straightforward to see that this mapping is a covering transformation and, furthermore, that it preserves the colour of the darts of $\Cov(G,\zeta)$. That is, each $f_g$ is an automorphism of the maniplex $\Cov(G,\zeta)$. This means that all the vertices within a fibre belong to the same orbit under the action of the (maniplex) automorphism group, and so $\Cov(G,\zeta)$ is an $n$-maniplex with at most $2$ orbits. If no automorphism of $\Cov(G,\zeta)$ preserving the colouring (that is, automorphisms of the resulting maniplex) maps a vertex in $\fib_u$ to a vertex in $\fib_v$, then $\Cov(G,\zeta)$ is a $2$-orbit $n$-maniplex and $\Gamma$ is its symmetry type graph. The non-existence of such an automorphism, however, will be somewhat intricate to prove.

The group $\overline{G}$ of automorphisms of $\Cov(G,\zeta)$ induced by right multiplication by elements of $G$ has  index at most $2$ in the full maniplex automorphism group $\Aut(\Cov(G,\zeta))$. Hence the orbits (on vertices) of $\overline{G}$, which correspond precisely to the fibres of $u$ and $v$, are blocks of imprimitivity. It follows that, if an automorphism $\tau$ maps one vertex in $\fib_u$ into a vertex in $\fib_v$, then it must map every vertex in $\fib_u$ to a vertex in $\fib_v$. This is, $\tau$ must be a lift of $R$, the unique colour preserving automorphism of $\Gamma$ that interchanges $u$ and $v$. Thus, to show that $\Cov(G,\zeta)$ is not a regular maniplex it suffices to show that $R$ does not lift. In light of the results shown in \cite{MNS}, we have that $R$ has a lift if and only $R$ maps isomorphically the group of closed walks based at $u$ into the group of closed walks based at $v$. The latter, in turn, implies that there is a induced group isomorphism $R^{\#}$ mapping the local voltage group $\Loc(G,u)$ into $\Loc(G,v)$ and satisfying $R^{\#}(\zeta(W)) = \zeta(R(W))$, for every closed walk $W$ based at $u$. Observe that, in the case of the pregraph $2_{\mI}^n$ with a voltage assignment like in Propositions \ref{prop:manilift1} or \ref{prop:manilift2}, we have $G \cong \Loc(G,u) \cong \Loc(G,v)$. Hence, we have the following two propositions 

\begin{proposition}
\label{prop:notau1}
Let $\Gamma$ be the pregraph $2_{\mI}^n$ with $1 \notin \mI$ and let $\zeta:D(\Gamma) \to G$ be the voltage assignment described in  Proposition \ref{prop:manilift1}. Let $R$ be the colour preserving automorphism interchanging the two vertices of $\Gamma$ .
If no automorphism of $G$ inverts $z$ while fixing $y_i$ for $i \in \{2,...,n-1\}$, then $R$ does not lift into an automorphism of $\Cov(\Gamma, \zeta)$.    
\end{proposition}

\begin{proposition}
\label{prop:notau2}
Let $\Gamma$ be the pregraph $2_{\mI}^n$ with $1 \in \mI$ and let $\zeta:D(\Gamma) \to G$ be the voltage assignment described in  Proposition \ref{prop:manilift2}. Let $R$ be the colour preserving automorphism interchanging the two vertices of $\Gamma$ .
If no automorphism of $G$ interchanges $z$ and $z'$ while fixing $y_i$ for $i \in \{2,...,n-1\}$, then $R$ does not lift into an automorphism of $\Cov(\Gamma, \zeta)$.        
\end{proposition}

\section{Construction}
\label{sec:last}

In this section we construct $2$-orbit $(n+1)$-maniplexes from regular $n$-maniplexes. In Theorems \ref{theo:1notinI} and \ref{theo:1inI} we show that for any $\mI \subsetneq \{1,\dots, n\}$ there is a maniplex of type $2^{n+1}_\mI$. This, together with Corollary \ref{cor:edges0n} proves Theorem \ref{theo:main}.

Let $\Gamma$ be the pregraph $2_{\mI}^{n+1}$ and let $\mathcal{M}$ be a regular $n$-polytope admitting a bi-colouring (white
and black) of flags consistent with $\mI$ (that is, $\Phi$ and $\Phi^{i}$ have the same colour if and only if $i\in \mI$). Let $\Mon_\mI(\mM)$ be the subgroup of $\Mon(\mM)$ preserving the colours of the bi-colouring and assume there exists a facet $F_0$ and $\eta\in \Mon_\mI(\mM)$ such 
that $\eta^2 = 1$ and
that for every two flags
$\Phi$ and $\Psi$ in $F$, $\Phi^\eta$ and $\Psi^\eta$
are in distinct facets of $\mM$. 
Now, choose an arbitrary white flag $\Phi_0$ in $F_0$ and define it to be the base flag of $F_0$.
Let $\Phi_1 = \Phi_0^\eta$ be the base flag of the facet $F_1$ containing $\Phi_1$.
Let $X=\mF(\mM) \times \ZZ_{2k}$ and $X_w$ the subset of $X$ consisting of pairs where the first element is a white flag.

We will provide a construction whose output is a permutation group acting on the set $X_w$. This group will serve as a voltage group to lift $\Gamma$ into a $(n+1)$-maniplex of type $2^{n+1}_\mI$.

The following two Remarks are a reminder of some of the contents of Section \ref{sec:intro}.

\begin{remark}\label{r:commute}
If $\mM$ is an $n$-maniplex and $r_0,r_1,\ldots,r_{n-1}$ are the standard generators of $\Mon(\mM)$ then
\begin{eqnarray*}
r_i^2 &=& 1 \cr
[r_i, r_j] &=& 1 \mbox{ if } |i-j|>1.
\end{eqnarray*}
\end{remark}

\begin{remark}
\label{rem:r0}
If $\Phi$ is a flag of a maniplex $\mM$, then $\Phi$ and $\Phi^{r_0}$ belong to the same facet of $\mM$.
\end{remark}

\begin{lemma}
\label{lem:distinctfacets}
With the notation above, the following hold:
\begin{enumerate}
\item If $\Phi\in F_0$, then $\Phi^\eta \not \in F_0$;
\item If $\Phi$ and $\Psi$ are distinct flags contained in the same facet of $\mM$, then
$\Phi^{\eta}$ and $\Psi^{\eta}$ are in distinct facets of $\mM$;

\item If $\Phi$ and $\Psi$ are distinct white flags contained in the same facet of 
$\mM$, then $\Phi^{r_0\eta r_0}$ and $\Psi^{r_0\eta r_0}$
         are in distinct facets of $\mM$.
\item If $\Phi$ is a white flag in $F_0$ other than $\Phi_0$, then $\Phi^{r_0\eta r_0}$ is neither in $F_0$ nor in $F_1$.
\end{enumerate}
\end{lemma}

\begin{proof}
In order to prove part (1), suppose that $\Phi^\eta \in F_0$ for some $\Phi\in F_0$ and let
$\Psi$ be an arbitrary flag in $F_0$ other than $\Phi$. Since $\mM$ is 
a regular polytope, there exists an automorphism $\gamma$ of $\mM$ mapping $\Phi$ to $\Psi$.
Note that $\gamma$ then preserves the facet $F_0$.
Since the actions of $\Mon(\mM)$ and $\Aut(\mM)$ commute, it follows that
 $\Psi^\eta = (\Phi^\gamma)^\eta = (\Phi^\eta)^\gamma$. Since $\Phi^\eta \in F_0$ and since
 $\gamma$ preserves $F_0$, it follows that $\Psi^\eta \in F_0$. This contradicts the assumptions on $\eta$.
 
Observe that part (2) follows directly from the regularity of $\mM$.

 To prove part (3),  observe that by Remark~\ref{rem:r0} and part (2), it follows that
 $\Phi^{r_0\eta}$ and $\Psi^{r_0\eta}$ are in distinct facets, and again by Remark~\ref{rem:r0} the result follows.

Finally, to prove part (4), observe that by part (1) and Remark~\ref{rem:r0}, $\Phi^{r_0\eta r_0}\not \in F_0$.
Further, since $\Phi^{r_0}$ is a black flag in $F_0$, it is not equal to $\Phi_0$, and by the definition of $\eta$,
$\Phi_0^\eta$ and  $\Phi^{r_0\eta}$ are in different facets of $\mM$. In particular, $\Phi^{r_0\eta}$ is not in $F_1$,
and by Remark~\ref{rem:r0} neither is $\Phi^{r_0\eta r_0}$.
\end{proof}

We will now fix a base flag in each facet of $\mM$. For an arbitrary white flag $\Phi$ of $F_0$, let $\Phi^{r_0\eta r_0}$ be the base flag of the facet that contains it.
Observe that all base flags determined so far are white, and that by Lemma~\ref{lem:distinctfacets} no facet is given more than one base flag, in this way.
For each remaining facet, choose its base flag to be an arbitrary white flag.

Let $\{r_{0},r_{1},...,r_{n-1}\}$ be the distinguished generators
of $\Mon(\mathcal{M})$. We extend the action of $r_i$ to $X$ by letting $(\Phi,\ell)^{r_i}=(\Phi^{r_i},\ell)$. 

For each facet $F$, let $\rho_0^F$ be the
automorphism of $F$ sending its distinguished base flag into its $0$-neighbour, and define $\widetilde{\rho_0}$ as the permutation on $\mF(\mM)$ that acts like $\rho_0^F$ on
the flags of each facet $F$. Note that the actions of $\trho$ and $r_0$ commute within each facet of $\mM$. We extend the action of $\trho$ to $X$ by letting $(\Phi,\ell)^{\trho} = (\Phi^{\trho},\ell)$ and therefore $r_0$ and $\trho$ commute as permutations of $X$.

Now let $y_{n}$ be the permutation of $X$ satisfying:
\begin{equation}
\label{eq:yn}
(\Phi, \ell)^{y_n}  =  \left\{ 
\begin{array}{ll}
(\Phi^{\trho r_0}, \ell+1) & \mbox{ if } \ell \mbox{ odd and } \Phi\in F_{0}, \mbox{ or } \ell \mbox{ even and } \Phi\in F_1; \\
(\Phi^{\trho r_0}, \ell-1) & \mbox{ if } \ell \mbox{ even and } \Phi\in F_{0}, \mbox{ or } \ell \mbox{ odd and } \Phi\in F_1;\\
(\Phi^{\trho r_0}, \ell) & \mbox{ otherwise.}
\end{array}
\right.
\end{equation}
Since $r_0$ and $\trho$ map white flags to black flags, $y_n$ permutes $X_w$.

In the following paragraphs we will define a set of generators for a voltage group $G \subset \mathcal{S}_{X_w}$ for $\Gamma$  and then prove that it has some desirable properties. Namely, that it agrees with the hypothesis of Propositions \ref{prop:manilift1} and \ref{prop:notau1} when $1\notin\mI$, or those of Propositions \ref{prop:manilift2} and \ref{prop:notau2} when $1 \in \mI$. We will produce two slightly different constructions yielding two different groups: one when $1 \in \mI$ and another when $1 \notin \mI$.
 
\subsection{Suppose $1 \notin \mI$}
 Define:
\begin{eqnarray*}
z&=&r_{0}r_{1},\\
y_{i}&=&\begin{cases}
r_{i} & \mbox{ if } i\in \mI,\\
r_{i}r_{0} & \mbox{ if } i \notin \mI.
\end{cases}
\mbox{ for } i \in \{2, \ldots, n-1\}, 
\end{eqnarray*}

Note that $z$ and $y_i$ are in $\Mon_\mI(\mM)$ and hence can be seen as permutations of $X_w$.

We define $G=\langle z,y_2,\dots,y_n \rangle$ and note that it acts on $X_w$.

\begin{lemma}\label{l:relations}
The group $G$ just defined together with its generators $z, y_2, \dots, y_n$ satisfy relations (\ref{eq:volt1}), (\ref{eq:volt2}) and (\ref{eq:volt3}).
\end{lemma}

\begin{proof}
For $i\in \{2,\ldots, n-1\}$, the relations (\ref{eq:volt1}) and (\ref{eq:volt3}) follow from Remark \ref{r:commute}.
To prove (\ref{eq:volt1}) for $i=n$, observe that $y_n^2=1$ follows directly from the definition of $y_n$. On the other hand, if $i<n-1$ then $\Phi$ and $\Phi^{r_i}$ belong to the same facet for every flag $\Phi$, and (\ref{eq:volt3}) follows from Remark \ref{r:commute} and the fact that $\trho$ commutes with $r_i$ for all $i \in \{0,\dots,n-2\}$.

Let us now prove relation (\ref{eq:volt2}). Suppose first that $i\in \{3,\ldots, n-1\}$ and recall that then $r_i$ commutes with $r_0$ and $r_1$.
If $i\in\mI$, then 
$z^{y_i} = r_ir_0r_1r_i =r_0r_1r_i^2 = z$.
If $i\not\in\mI$, then 
$z^{y_i} = r_ir_0^2r_1r_ir_0 =r_1r_0r_i^2 = z^{-1}$.
In order to prove the relation (\ref{eq:volt2}) for $i=n$, recall that for any flag $\Phi$
the flags $\Phi^z$ and $\Phi$ belong to the same facet. The fact that $z^{y_n} = z^{-1}$  follows directly from the definition of $y_n$ and
the fact that $r_1r_0 = z^{-1}$.
\end{proof}

In what follows we will define an element $\mu$ of $G$ and we will show that, if an automorphism $\tau \in \Aut(G)$ fixes all $y_i$'s while inverting $z$, then $\mu$ and its image under $\tau$ cannot have the same order. Thus proving that such an automorphism cannot exists. 

\begin{figure}[h!]
\centering
\includegraphics[width=1\textwidth]{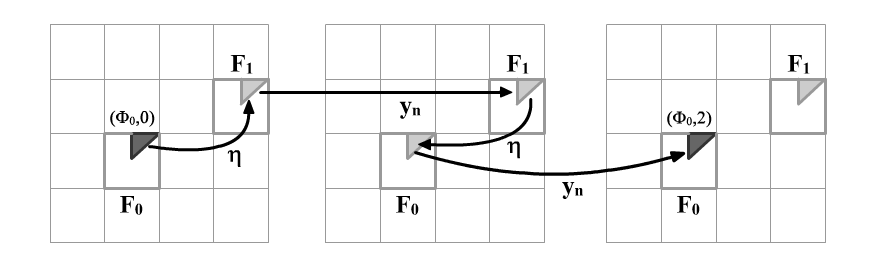}
\caption{}
\label{fig:mu}
 \end{figure}

 Let $\mu=(\eta y_{n})^2$. Now let $\ell$ be an arbitrary odd element of $\ZZ_{2k}$.
Then:

$$
 (\Phi_0, \ell)^\mu=(\Phi_0^\eta, \ell)^{y_n\eta y_n}=  (\Phi_1, \ell)^{y_n\eta y_n}= (\Phi_1^{r_0\trho}, \ell-1)^{\eta y_n} 
$$ 
Since $\Phi_1$ is the base flag of the facet $F_1$, the permutation $r_0\trho$ fixes it. Therefore,
the above equals:
 $$
  (\Phi_1, \ell-1)^{\eta y_n} =
 (\Phi_1^\eta, \ell-1)^{ y_n}=
 (\Phi_0,\ell-1)^{y_n}=  (\Phi_0^{r_0\trho},\ell-2) = (\Phi_0,\ell-2),
 $$
as illustrated in Figure \ref{fig:mu}.
By applying this argument $k$ times, we get the following lemma. The orbit of $\langle \mu \rangle$ containing $(\Phi_0, 1)$
consists of all the elements of the form $(\Phi_0, \ell)$ for $\ell$ odd and thus:

\begin{lemma}
\label{lem:orbitsize}
The orbit of $\langle \mu \rangle$ containing $(\Phi_0, 1)$ has size $k$.
\end{lemma}

We want to
show that no automorphism of $G=\langle z, y_{2},y_{3,...,}y_{n-1},y_{n} \rangle$
fixes all generators $y_{i}$ while inverting $z$. Suppose one such
automorphism, $\tau$, exists. Observe that the action of $\tau$ upon the
subgroup $H=\langle y_{2},y_{3,...,}y_{n-1},z \rangle \le G$ coincides with the conjugation by $r_0$.
In particular, since $\eta \in H$ and $y_n^\tau = y_n$, we see that
\begin{equation}
\label{eq:mutau}
\mu^\tau = ((\eta y_n)^2)^\tau = (\eta^{r_0}y_n)^2 = (r_0 \eta r_0 y_n)^2 .
\end{equation}

For $\ell\in \ZZ_{2k}$ let $X_\ell = \{ (\Phi,\ell) : \Phi \in F(\mM)\}$.
The aim of the following three lemmas is to show that an orbit
of every element $(\Phi,\ell) \in X_w$ under  $\langle \mu^\tau\rangle$
is contained either in $X_\ell \cup X_{\ell-1}$ or in $X_\ell \cup X_{\ell+1}$.
This will then imply that for a sufficiently large $k$, 
the order of $\mu$ is larger than the order of $\mu^\tau$, contradicting our assumption on $\tau$.

\begin{lemma}
\label{lem:r0trho}
Let $\Phi$ be a white flag in $F_0$. Then
$\Phi^\eta$ and  $\Phi^{r_0 \eta r_0}$ are fixed by $\trho r_0$.
\end{lemma}

\begin{proof}
The flags $\Phi^\eta$ and $\Phi^{r_0 \eta}$ are base flags of some facets of $\mM$, and therefore $\trho$ acts like $r_0$ on each of them. Then $\Phi^{\eta \trho r_0} =\Phi^{\eta}$, and since $\trho$ and $r_0$ commute, $(\Phi^{r_0 \eta r_0})^{\trho r_0} = \Phi^{r_0 \eta \trho} = \Phi^{r_0 \eta r_0}$.
\end{proof}

\begin{lemma}
\label{lem:cases}
Let $\Phi$ be a white flag, let $\ell \in \ZZ_{2k}$. Let $(\Psi_1,\ell_1) = (\Phi,\ell)^{r_0\eta r_0 y_n r_0 \eta r_0}$ and
let $(\Psi_2,\ell_2) = (\Phi,\ell)^{r_0\eta r_0}$. If $\Psi_1$ or $\Psi_2$ is in $F_0$,
then the orbit of $(\Phi,\ell)$ under $\langle \mu^\tau \rangle$ is of length $2$.

\end{lemma}

\begin{proof}
Suppose first that $\Psi_1$ is in $F_0$. Let $(\Psi,\ell') =  (\Psi_1,\ell_1)^{y_n}$ 
and observe that $\Psi \in F_0$ and that $(\Psi,\ell') = (\Phi,\ell)^{\mu^\tau}$.
We will now show that the orbit of $(\Psi,\ell')$ under $\langle \mu^\tau \rangle$ has length $2$, implying that so does the orbit of
$(\Phi,\ell)$.

By part (4) of Lemma~\ref{lem:distinctfacets}, 
it follows that $\Psi^{r_0 \eta r_0}$ is neither in $F_0$ nor in $F_1$. 
Therefore $(\Psi^{r_0 \eta r_0},\ell')^{y_n} =(\Psi^{r_0 \eta r_0 \trho r_0},\ell')$, and by Lemma \ref{lem:r0trho}, this equals $(\Psi^{r_0 \eta r_0},\ell')$. Hence
$$
(\Psi,\ell')^{\mu^\tau}  = (\Psi,\ell')^{(r_0 \eta r_0y_n)^2} =(\Psi^{r_0 \eta r_0},\ell')^{y_nr_0\eta r_0y_n} =
(\Psi^{r_0 \eta r_0},\ell')^{r_0\eta r_0y_n} = (\Psi,\ell')^{y_n}
$$
By the definition of $y_n$ the latter equals $(\Psi,\ell' +1)$ if $\ell'$ is odd, and $(\Psi,\ell' -1)$ if $\ell'$ is even.
By applying $\mu^\tau$ to $(\Psi,\ell')^{\mu^\tau}$, we thus obtain $(\Psi,\ell)$.
 This shows that the orbit $(\Psi,\ell')$ under the action of $\langle \mu^\tau \rangle$ contains precisely two elements.
\medskip

Suppose now that $\Psi_2$ is in $F_0$ and recall that $(\Psi_2,\ell_2) = (\Phi^{r_0\eta r_0},\ell)$.
 Then $(\Phi^{r_0 \eta r_0})^{r_0 \eta r_0} = \Phi$ is the base flag of
some facet of $\mM$, and by part (4) of Lemma~\ref{lem:distinctfacets}, $\Phi$ belongs to neither $F_0$ nor $F_1$.

Therefore 
\begin{equation}
\label{eq:1}
(\Phi^{r_0 \eta r_0},\ell)^{y_n} =(\Phi^{r_0 \eta r_0 (\trho r_0)},\ell+(-1)^\epsilon)
\end{equation}
where $\epsilon$ is $0$ or $1$ depending
on whether  $\ell$ is odd or even, respectively.
Then
$$
(\Phi,\ell)^{\mu^\tau}  = (\Phi,\ell)^{(r_0 \eta r_0y_n)^2} =(\Phi^{r_0 \eta r_0},\ell)^{y_nr_0\eta r_0y_n} =
(\Phi^{r_0 \eta r_0(\trho r_0)(r_0\eta r_0)},\ell+(-1)^\epsilon)^{y_n}
$$
Since $\Phi^{r_0 \eta r_0} \in F_0$ and since $\trho r_0$ preserves the set of flags of $F_0$, we see that
$\Phi^{r_0 \eta r_0(\trho r_0)} \in F_0$ and hence $\Phi^{r_0 \eta r_0(\trho r_0)(r_0\eta r_0)}$ is the base flag
of some facet in $\mM$, which by part (4) of Lemma~\ref{lem:distinctfacets} is neither $F_0$ nor $F_1$.
By the definition of $y_n$, Lemma~\ref{lem:r0trho} and the computation above, it follows that
$$
(\Phi,\ell)^{\mu^\tau}  = (\Phi^{r_0 \eta r_0(\trho r_0)(r_0\eta r_0)},\ell+(-1)^\epsilon).
$$
 Let us now apply $\mu^\tau$ one more time. Note that in the computation below we will make use of formula~(\ref{eq:1})
 and the fact that $y_n$ is an involution.
$$
(\Phi,\ell)^{(\mu^\tau)^2}  = (\Phi^{r_0 \eta r_0(\trho r_0)(r_0\eta r_0)},\ell+(-1)^\epsilon)^{\mu^\tau} = 
(\Phi^{r_0 \eta r_0(\trho r_0)(r_0\eta r_0)},\ell+(-1)^\epsilon)^ {(r_0 \eta r_0y_n)^2}
$$
$$
 = (\Phi^{r_0 \eta r_0(\trho r_0)},\ell+(-1)^\epsilon)^ {y_nr_0 \eta r_0y_n} =
(\Phi^{r_0 \eta r_0},\ell)^ {r_0 \eta r_0y_n} = (\Phi,\ell)^ {y_n}.
$$
Recall that $\Phi$ is the base flag of a facet
that is neither $F_0$ nor $F_1$. The definition of $y_n$ and
Lemma~\ref{lem:r0trho} then imply that $(\Phi,\ell)^{y_n} = (\Phi,\ell)$. In particular, 
the orbit of $(\Phi,\ell)$ under $\langle \mu^\tau\rangle$ again consists of two elements only.

\end{proof}

\begin{lemma}
\label{lem:new}
For an element $\ell \in \ZZ_{2k}$ let $t$ be such that $\ell \in \{2t,2t+1\}$.
Then:
\begin{enumerate}
\item
 If $(\Phi,\ell)^{\mu^\tau} \not \in X_{2t}\cup X_{2t+1}$,
then either $\Phi^{r_0\eta r_0}$ or $\Phi^{r_0\eta r_0 y_n r_0 \eta r_0}$ is in $F_0$.
\item
 If the orbit of $(\Phi,\ell)$ under $\langle \mu^\tau \rangle$ is not contained in $X_{2t}\cup X_{2t+1}$,
then its length is $2$.
\end{enumerate}
\end{lemma}

\begin{proof}
Part (1) follows directly from formula (\ref{eq:mutau}) and the definition of $y_n$.
To prove part (2), suppose that the orbit of $(\Phi,\ell)$ under $\langle \mu^\tau \rangle$ is of length at least $3$ and is
not contained in $X_{2t}\cup X_{2t+1}$.
Let $i$ be the smallest integer such that $(\Phi,\ell)^{(\mu^\tau)^{i+1}} \not \in X_{2t}\cup X_{2t+1}$.
Let $(\Psi,\ell') = (\Phi,\ell)^{(\mu^\tau)^{i}}$. Then $\ell' \in \{2t,2t+1\}$ and by
 part (1), either $\Psi^{r_0\eta r_0}$ or $\Psi^{r_0\eta r_0 y_n r_0 \eta r_0}$ is in $F_0$. By Lemma~\ref{lem:cases},
 the orbit of $(\Psi,\ell')$ under $\langle \mu^\tau \rangle$ is of length $2$, and thus so is the orbit of $(\Phi,\ell)$.
\end{proof}

\begin{theorem}
\label{theo:1notinI}
The group $G=\langle z, y_2, \dots, y_n \rangle$ defined as above satisfies relations (\ref{eq:volt1}), (\ref{eq:volt2}) and (\ref{eq:volt3}), and for a sufficiently large $k$ (in the definition of $X$) it admits no automorphism that inverts $z$ while fixing $y_i$ for $i \in \{2, \dots, n\}$. 
\end{theorem}

\begin{proof}
Lemma \ref{l:relations} shows that $G$ satisfies the desired relations.
Now suppose there exists an automorphism $\tau$ of $G$ that inverts $z$ while fixing $y_i$ for $i \in \{2, \dots, n\}$. It follows from Lemma \ref{lem:new} that the orbit of any pair $(\Phi,\ell)$ under $\langle \mu^\tau \rangle$ is of size at most $2|\mF_w(\mM)|$, where $\mF_w(\mM)$ denotes the set of white flags of $\mM$. Recall that the orbit of $(\Phi_0,1)$ under $\langle \mu \rangle$ is of size $k$. By choosing $k>2|\mF_w(\mM)|$ in the definition of $X$, we ensure that order of $\mu$ is necessarily larger than the order of $\mu^\tau$, contradicting that $\tau$ is an automorphism of $G$. The result follows.

\end{proof}

\subsection{Suppose $1 \in \mI$}
 Define:

\begin{eqnarray*}
z&=&r_{1},\\
z'&=&r_{0}r_{1}r_{0},\\
y_{i}&=&\begin{cases}
r_{i} & \mbox{ if } i\in \mI,\\
r_{i}r_{0} & \mbox{ if } i \notin \mI.
\end{cases}
\mbox{ for } i \in \{2, \ldots, n-1\}, 
\end{eqnarray*}

Note that $z$, $z'$ and $y_i$ are in $\Mon_\mI(\mM)$ and hence can be seen as permutations of $X_w$. 
Define $G=\langle z,z',y_2,\dots,y_n \rangle$ and note that it acts on $X_w$.

\begin{lemma}\label{lem:relations2}
The group $G$ just defined together with its generators $z,z', y_2, \dots, y_n$ satisfy relations (\ref{eq:1inivolt1}), (\ref{eq:1inivolt2}), (\ref{eq:1inivolt3}) and (\ref{eq:1inivolt4}).
\end{lemma}

\begin{proof}
From remark \ref{r:commute} we have that $z^2=(z')^2=y_i=1$, when $i < n$. That $y_n^2=1$ follows from the definition. Hence, condition (\ref{eq:1inivolt1}) holds.

To show that condition (\ref{eq:1inivolt2}) is satisfied first suppose $i \in \mI$. Then $z^{y_i}=y_izy_i=r_ir_1r_i=r_1=z$. Now, if $i \notin \mI$ and $i<n$, then $z^{y_i}=r_ir_0r_1r_ir_0=r_ir_ir_0r_1r_0=r_0r_1r_0=z'$. Finally, note that $\trho$ commutes with $r_0$ and $r_1$. Thus $z^{y_n}=z'$. 

Let us now prove relation (\ref{eq:1inivolt3}). Suppose $i \in \mI$, so $(z')^{y_i}=r_ir_0r_1r_0r_i=r_ir_ir_0r_1r_0=z'$. If $i \notin \mI$ and $i<n$, this is just equivalent to condition (\ref{eq:1inivolt2}) when $i<n$. Once more, since $\trho$ commutes with each $r_i$, $i \in \{0,1\}$, we have that $(z')^{y_n}=z$. Condition (\ref{eq:1inivolt3}) is therefore satisfied.

Now, if $i,j<n$, from Remark \ref{r:commute} we have that $[y_i,y_j]=1$. On the other hand, if $i=n$, condition (\ref{eq:1inivolt4}) follows from the fact that $\trho$ commutes with $r_i$ for all $i \in \{0,...,n-2\}$.
\end{proof}

 Let $\mu=(\eta y_{n})^2$, just as in the previous case.

We want to show that no automorphism of $G=\langle z,z', y_{2},y_{3,...,}y_{n-1},y_{n} \rangle$
fixes all generators $y_{i}$ while interchanging $z$ and $z'$. Suppose one such
automorphism, $\tau$, exists. Again, the action of $\tau$ upon the
subgroup $H=\langle y_{2},y_{3,...,}y_{n-1},z,z' \rangle \le G$ coincides with the conjugation by $r_0$. Notice that Lemmas \ref{lem:r0trho}, \ref{lem:cases} and \ref{lem:new} do not depend on $z$ nor $z'$, but only on the definitions of $y_n$, $\eta$ and the fact that $\tau$ acts like conjugation by $r_0$. Therefore, Lemmas \ref{lem:r0trho}, \ref{lem:cases} and \ref{lem:new} also hold when $1 \in \mI$. We thus have the following theorem.

\begin{theorem}
\label{theo:1inI}
The group $G=\langle z, z', y_2, \dots, y_n \rangle$ defined as above satisfies relations (\ref{eq:1inivolt1}), (\ref{eq:1inivolt2}), (\ref{eq:1inivolt3}) and (\ref{eq:1inivolt4}), and for a sufficiently large $k$ (in the definition of $X$) it admits no automorphism that interchanges $z$ and $z'$ while fixing $y_i$ for $i \in \{2, \dots, n\}$. 
\end{theorem}

\begin{proof}
Lemma \ref{lem:relations2} shows that $G$ satisfies the desired relations.
The proof that no automorphism of $G$ fixes $y_i$, $i \in \{2,...,n\}$, while interchanging $z$ and $z'$ is identical to that of Theorem \ref{theo:1notinI}.

\end{proof}

We are finally ready to prove Theorem \ref{theo:main}. The proof is by induction over $n$, where the base case $n=3$ is known to be true (see for example \cite{Isa10}).

Assume that the result is true for certain $n$ and let $\Gamma$ be the pregraph $2_{\mI}^{n+1}$. The result follows from Corollary \ref{cor:edges0n} if $0 \in \mI$, and from Theorems \ref{theo:1notinI} and \ref{theo:1inI} if $0 \notin \mI$.

\section{Acknowledgements}
The first author was partially supported by PAPIIT-UNAM under grant IN-107015 and
the second and third author by the Slovenian Research Agency (the {\em Young researcher scholarship} and the 
research core funding P1-0294, respectively).

\end{document}